\documentclass[10pt]{amsart}
\usepackage{latexsym}
\usepackage{amsmath,amsthm,amsfonts,amscd,eucal}
\usepackage{epsfig}

\newtheorem{thm}{Theorem}[section]
\newtheorem{cor}[thm]{Corollary}
\newtheorem{lem}[thm]{Lemma}
\newtheorem{prop}[thm]{Proposition}

\theoremstyle{definition}
\newtheorem{defn}[thm]{Definition}
\theoremstyle{remark}
\newtheorem{rem}[thm]{Remark}
\newtheorem{ex}[thm]{Example}
\numberwithin{equation}{section}

\def\ca{{\mathcal A}}
\def\cb{{\mathcal B}}

\def\ce{{\mathcal E}}
\def\cf{{\mathcal F}}
\def\cg{{\mathcal G}}

\def\cp{{\mathcal P}}

\def\bc{{\mathbb C}}

\newcommand{\bn}{\mathbb N}
\def\br{{\mathbb R}}
\def\bz{{\mathbb Z}}

\renewcommand{\a}{\alpha}
\def\b{\beta}
\def\g{\gamma}        \def\G{\Gamma}
        \def\D{\Delta}
\def\eps{\varepsilon}

\def\th{\vartheta}

\def\l{\lambda}       
\def\m{\mu}
\def\n{\nu}

\def\r{\rho}
\def\s{\sigma}       
    
\def\t{\tau}

\def\f{\varphi}

        \def\O{\Omega}

\newcommand{\norm}[1]{\left\Vert#1\right\Vert}

\newcommand{\set}[1]{\left\{#1\right\}}
\DeclareMathOperator{\diam}{diam}

\newcommand{\itm}[1]{\item[($#1$)]}
\newcommand{\de}{\partial}

\newcommand{\supp}{\text{supp}}
\newcommand{\cell}{\ce}
\newcommand{\Tr}{Tr^{\cg}}

\def\ov{\overline}

\begin{document}

\title[Self-similar CW-complexes and L$^{2}$-invariants]
{A C$^{*}$-algebra of geometric operators \\ on self-similar 
CW-complexes.\\ Novikov-Shubin and L$^{2}$-Betti numbers}
\author{Fabio Cipriani, Daniele Guido, Tommaso Isola}%
\date{\today}
\address{(F.C.) Dipartimento di Matematica, Politecnico di Milano,
piazza Leonardo da Vinci 32, 20133 Milano, Italy. \\
 (D.G., T.I.) Dipartimento di Matematica, Universit\`a di Roma ``Tor
Vergata'', I--00133 Roma, Italy.}%
\email{fabcip@mate.polimi.it, guido@mat.uniroma2.it,
isola@mat.uniroma2.it}%
\thanks{This work has been partially supported by GNAMPA, MIUR and by
the European Networks ``Quantum Spaces - Noncommutative Geometry"
HPRN-CT-2002-00280, and ``Quantum Probability and Applications to
Physics, Information and Biology''}

\begin{abstract}
   A class of CW-complexes, called self-similar complexes, is
   introduced, together with C$^{*}$-algebras $\ca_{j}$ of operators,
   endowed with a finite trace, acting on square-summable cellular
   $j$-chains.  Since the Laplacian $\D_{j}$ belongs to $\ca_{j}$,
   L$^2$-Betti numbers and Novikov-Shubin numbers are defined for
   such complexes in terms of the trace.  In particular a relation involving the
   Euler-Poincar\'e characteristic is proved.  L$^2$-Betti and
   Novikov-Shubin numbers are computed for some self-similar complexes
   arising from self-similar fractals.
\end{abstract}
\maketitle

\section{Introduction.}

In this paper we address the question of the possibility of extending
the definition of some $L^{2}$-invariants, like the L$^2$-Betti
numbers and Novikov-Shubin numbers, to geometric structures which are
not coverings of compact spaces.

The first attempt in this sense is due to John Roe \cite{RoeIndex},
who defined a trace on finite propagation operators on amenable
manifolds, allowing the definition of L$^2$-Betti numbers on these
spaces.  However such trace was defined in terms of a suitable
generalised limit, hence the corresponding L$^2$-Betti and
Novikov-Shubin numbers also depend on this generalised limit
procedure.

Here we show that, on spaces possessing a suitable self-similarity, 
it is possible to select a natural C$^{*}$-algebra of operators, generated by 
operators with finite propagation and locally commuting with the 
transformations giving the self-similar structure, on which a Roe-type
trace is well defined. 


The theory of $L^{2}$-invariants was started by Atiyah, who, in a
celebrated paper \cite{Atiyah}, observed that on covering manifolds
$\G\to M \to X$, a trace on $\Gamma$-periodic operators may be
defined, called $\Gamma$-trace, with respect to which the Laplace
operator has compact resolvent.  Replacing the usual trace with the
$\Gamma$-trace, he defined the $L^2$-Betti numbers and proved an index
theorem for covering manifolds.

Based on this paper, Novikov and Shubin \cite{NS1} observed that,
since for noncompact manifolds the spectrum of the Laplacian is not
discrete, new global spectral invariants can be defined, which
necessarily involve the density near zero of the spectrum.
 
L$^{2}$-Betti numbers were proved to be $\G$-homotopy invariants by
Dodziuk \cite{Dodziuk}, whereas Novikov-Shubin numbers were
proved to be $\G$-homotopy invariants by Gromov-Shubin \cite{GS}. 
L$^{2}$-Betti numbers (depending on a generalised limit procedure) 
were subsequently defined for open manifolds by
Roe, and were proved to be invariant under quasi-isometries
\cite{RoeBetti}. The invariance of Novikov-Shubin numbers was proved 
in \cite{GuIs4}.


The basic idea of the present analysis is the notion of self-similar
CW-complex, which is defined as a complex endowed with a natural
exhaustion $\set{K_{n}}$ in such a way that $K_{n+1}$ is a union (with
small intersections) of a finite number of copies of $K_{n}$.  The
identification of the different copies of $K_{n}$ in $K_{n+1}$ gives
rise to many local isomorphisms on such complexes.  Then we consider
finite propagation operators commuting with these local isomorphisms
up to boundary terms, and call them geometric operators.  Geometric operators generate a C$^{*}$-algebra
$\ca_j$ on the space of $\ell^{2}$-chains of $j$-cells, for any $j$ from zero to the dimension of the complex, containing the $j$-Laplace operator.  For any operator $T$ in this
C$^{*}$-algebra, we consider the sequence of the traces of $TE_{n}$,
renormalised with the volume of the $j$-cells of $K_{n}$, where
$E_{n}$ denotes the projection onto the space generated by  the $j$-cells of $K_{n}$.  Such a
sequence is convergent, and the corresponding functional is indeed a
finite trace on $\ca_j$.  By means of these traces, L$^2$-Betti
numbers and Novikov-Shubin numbers are defined. For the sake of completeness, we mention that notions related to that of geometric operators have been considered in the literature, see $e.g.$ \cite{LeSt} where they are called tight binding operators, and \cite{Elek1} where they are called pattern invariant operators.

In the $\Gamma$-covering case, $L^{2}$-Betti numbers are defined as
$\Gamma$-dimensions of the kernels of Laplace operators, namely as
$\Gamma$-traces of the corresponding projections.  This is not allowed
in our framework.  Indeed, our traces being finite, and the
C$^{*}$-algebras being weakly dense in the algebra of all bounded
linear operators, our traces cannot extend to the generated von
Neumann algebras.  In particular they are not defined on the spectral
projections of the Laplace operators.  Therefore we define L$^2$-Betti
numbers as the infimum of the traces of all continuous functional
calculi of the Laplacian, with functions taking value 1 at 0, namely
L$^2$-Betti numbers are defined as the ``external measure'' of the
spectral projections of the Laplace operators.

Since we are in an infinite setting, the Euler-Poincar\'e
characteristic is naturally defined as a renormalised limit of the
Euler-Poincar\'e characteristic of the truncations $K_{n}$ of the
complex.  We prove that such characteristic coincides with the
alternating sum of the L$^2$-Betti numbers.  An analogous result,
though obtained with a different proof, for amenable simplicial
complexes is contained \cite{Elek}.

Here we do not prove directly invariance results for L$^2$-Betti or
Novikov-Shubin numbers, however when 1-dimensional CW-complexes are
considered, and in particular prefractal graphs determined by nested
fractals, a result by Hambly and Kumagai applies \cite{HaKu}, implying
that Novikov-Shubin numbers are invariant under rough isometries.
Further results on invariance will be proved elsewhere \cite{CFGI}.


We then show that in some cases L$^2$-Betti and Novikov-Shubin numbers
can be computed, relying on results of several authors concerning
random walks on graphs.  In particular, it turns out that the
Novikov-Shubin numbers of some prefractal complexes coincide with the
spectral dimensions of the corresponding fractals, thus strenghtening
the interpretation of such numbers as (asymptotic) spectral dimensions
given in \cite{GuIs4}.


Our framework was strongly influenced by the approach of Lott and
L\"{u}ck \cite{LoLu}, in particular we also consider invariants
relative to the boundary, however we are not able to prove the
Poincar\'e duality shown in \cite{LoLu}.

The paper is organised as follows. In Section 2 we recall some 
notions from the theory of CW-complexes and introduce the basic 
operators. Section 3 introduces the notion of local isomorphisms of 
CW-complexes and the algebra of geometric operators. The notion of 
self-similar CW-complex is given in Section 4, and  a finite trace on geometric operators is constructed.

In Section 5 we introduce L$^{2}$-Betti and Novikov-Shubin numbers for
the above setting, and prove the mentioned result on the
Euler-Poincar\'e characteristic.  Section 6 focuses on the subclass of
self-similar CW-complexes given by prefractal complexes, and on some
properties of the associated Laplacians.  Computations of the
Novikov-Shubin numbers for fractal graphs in terms of transition
probabilities, together with an invariance result under rough
isometries are discussed in Sections 7 and 8, and the top-dimensional
relative Novikov-Shubin number is computed for two examples of
2-dimensional CW-complexes.

In closing this introduction, we note that the C$^*$-algebra and the trace for self-similar graphs constructed in this paper, are used in \cite{GILa01} to study the Ihara zeta function for fractal graphs.

The results contained in this paper were announced in the Conferences
``C$^{*}$-algebras and elliptic theory'' Bedlewo 2006, and
``21$^{st}$ International Conference on Operator Theory'' Timisoara 2006.

\section{CW-complexes and basic operators.}

 In this paper we shall consider a particular class of infinite
 CW-complexes, therefore we start by recalling some notions from
 algebraic topology, general references being \cite{LuWe,Munkres}.  A
 {\it CW-complex} $M$ of dimension $p\in \mathbb{N}$ is a Hausdorff
 space consisting of a disjoint union of (open) cells of dimension
 $j\in\set{0, 1, \dots , p}$ such that: $(i)$ for each $j$-cell
 $\s^{j}_{\a}$, there is a continuous map $f^{j}_{\a}: \{x\in\br^{j} :
 \|x\|\leq 1\} \to X$ that is a homeomorphism of $\{x\in\br^{j} :
 \|x\| < 1\}$ onto $\s^{j}_{\a}$, and maps $\{x\in\br^{j} : \|x\| =
 1\}$ into a finite union of cells of dimension $<j$; $(ii)$ a set
 $A\subset X$ is closed in $X$ iff $A\cap \ov{\s}^{j}_{\a}$ is closed
 in $\ov{\s}^{j}_{\a}$, for all $j$, $\a$, where $\ov{\s}^{j}_{\a}$
 denotes the closure of ${\s}^{j}_{\a}$ in $M$.  Let us denote by
 $\dot{\s}^{j}_{\a}= f^{j}_{\a}(\set{x\in\br^{j} : \|x\| =1})$ the
 boundary of $\s^{j}_{\a}$, for all $j$, $\a$.  A CW-complex is {\it
 regular} if $f^{j}_{\a}$ is a homeomorphism, for all $j$, $\a$.

 We denote by $\cell_{j}(M) := \{ \s^{j}_{\a} : \a\in\ca_{j} \}$,
 $j=0, 1, \dots , p$, the family of $j$-cells, and by $M^{j} :=
 \cup_{k=0}^{j} \cell_{k}(M)$, the $j$-skeleton of $M$.  Then $C_j
 (M):= H_{j}(M^{j},M^{j-1},\bz)$ is the (abelian) group of
 $j$-dimensional cellular chains, and is generated by the class of
 $\s^{j}_{\a}$, $\a\in\ca_{j}$.  Let $\de_{j}:C_{j}(M)\to C_{j-1}(M)$
 be the boundary operator, which is the connecting homomorphism of the
 homology sequence of the triple $(M^{j},M^{j-1},M^{j-2})$.  
 Let us choose an orientation of $M$, that is, a basis $\{
 \widehat{\s}^{j}_{\a}: \a\in\ca_{j} \}$ of $C_j (M)$,
 $j\in\set{0,\ldots,p}$, where each $\widehat{\s}^{j}_{\a}$ is (up to
 sign) the class of one (open) $j$-cell.  We will usually identify the
 algebraic cell $\widehat{\s}^{j}_{\a}$ with the geometric cell
 $\s^{j}_{\a}$, and denote by $-\s^{j}_{\a}$ the cell $\s^{j}_{\a}$
 with the opposite orientation.  Then the action of $\de_{j}$ on the
 chosen basis is given by $\displaystyle\de_{j} \s^{j}_{\a} =
 \sum_{\b\in\ca_{j-1}} [\s^{j}_{\a}:\s^{j-1}_{\b}] \s^{j-1}_{\b}$,
 where $[\s^{j}_{\a}:\s^{j-1}_{\b}]\in\bz$ depends on the chosen
 orientation and is called incidence number.  If $M$ is regular,
 $[\s^{j}_{\a}:\s^{j-1}_{\b}] \in \set{-1,0,1}$, and
 $[\s^{j}_{\a}:\s^{j-1}_{\b}] = 0 \iff \s^{j-1}_{\b}\cap
 \ov{\s}^{j}_{\a} =\emptyset$.
 Let us recall that the orientation of the zero-cells is chosen in
 such a way that, for any 1-cell $\s^{1}$,
 $\sum_{\a}[\s^{1},\s_{\a}^{0}]=0$.

 In the following we will consider only regular CW-complexes, unless
 otherwise stated.

 A Hilbert norm on $C_j (M)\otimes_{\bz} \bc$ is then defined as
 $\|c\|^{2}:= \sum_i |c_i|^2$ when $c=\sum_i c_i\cdot \sigma_i\in
 C_j(M)\otimes_{\bz} \bc$.  The Hilbert space $C^{(2)}_{j}(M) \equiv
 \ell^2 (\cell_j M)$ is the completion of $C_j (M)\otimes_{\bz}\bc$
 under this norm.

 \noindent We can extend $\de_{j}$ to a densely defined linear
 operator $C^{(2)}_{j}(M) \to C^{(2)}_{j-1}(M)$.  Then the {\it
 half--Laplace operators} $\Delta_{j\pm}$ are
\begin{equation*}
\begin{split}
\Delta_{j+}:&= \partial_{j+1} \partial_{j+1}^*\\
\Delta_{j-}:&=\partial_{j}^* \partial_{j}
\end{split}
\end{equation*}
 and the {\it Laplace operators} are $\Delta_j := \Delta_{j+} +
 \Delta_{j-}$.  These are operators on $\ell^2 (\cell_j M)$ densely
 defined on $C_j (M)\otimes_{\bz}\bc$.

 Let us observe that $\partial_j$ is a bounded operator under some
 condition.

 \begin{defn}[Bounded complex]  Let $M$ be a regular
 CW-complex, denote by
     \begin{align*}
     V_{j}^{+} & := \sup_{\s\in\cell_{j}(M)}
     |\set{\t\in\cell_{j+1}(M) : \dot{\t}\supset\s }| \\
     V_{j}^{-} & := \sup_{\s\in\cell_{j}(M)}
     |\set{\r\in\cell_{j-1}(M) : \r\subset\dot{\s} }| ,
     \end{align*}
 where $|\cdot|$ denotes the cardinality. We say that $M$ is a 
bounded complex if
 $V_{j}^{\pm}<\infty$, for all $j$.
 \end{defn}

 \begin{lem}\label{Lemma:boundedLaplacians}
     Let $M$ be a bounded regular CW-complex.  Then
     $\partial_j : \ell^2 (\cell_j M) \to \ell^2 (\cell_{j-1} M)$ is
     bounded.
\end{lem}
\begin{proof}
    If $c=\sum_i c_i\cdot \sigma_i$, setting $\alpha
    \sim\beta$ if there is $\r\in\cell_{j-1}(M)$ s.t. $\r\subset
    \dot{\s}_{\a}\cap\dot{\s}_{\b}$, we have
    \begin{equation*}
    \begin{split}
	\|\partial_j c\|^2 &= \sum_{\alpha \sim\beta} \ov{c_\alpha}
	\cdot c_\beta \cdot (\de_j\sigma_\alpha ,\de_j\sigma_\beta )\\
	&\le \sum_{\alpha \sim\beta} |c_\alpha | \cdot |c_\beta |\cdot
	|(\de_j\sigma_\alpha ,\de_j\sigma_\beta )|\\
	& \le \frac12 V_j^{-} \sum_{\alpha \sim\beta} \Bigl( |c_\alpha
	|^2 + |c_\beta |^2 \Bigr)\\
	&\le (V_j^{-})^2 V_{j-1}^{+} \|c\|^2 .
    \end{split}
    \end{equation*}
    Indeed 
    \begin{align*}
	|(\de_j\sigma_\alpha ,\de_j\sigma_\beta )| & \leq
	\sum_{\r\subset\dot{\s}_{\a}\cap\dot{\s}_{\b}} |[\s_{\a}:\r]|
	\cdot |[\s_{\b}:\r]| \\
	& \leq |\set{\r\in\cell_{j-1}:
	\r\subset\dot{\s}_{\a}\cap\dot{\s}_{\b}} | \leq V_{j}^{-},
    \end{align*}
    while, for any $\a\in\ca_{j}$, $\displaystyle |\set{ \b\in\ca_{j}
    : \b\sim\a }| \leq V_{j}^{-}V_{j-1}^{+}$.
\end{proof}

\begin{lem}
    $$
    \partial^{*}_{j+1}\sigma = \sum_{\t\in\cell_{j+1}(M)}
    [\t:\s]\tau\, .
    $$
\end{lem}
\begin{proof}
    Indeed, with $\t\in\cell_{j+1}(M)$,
    \[
    (\t, \de^{*}_{j+1} \s) = (\de_{j+1} \t,\s) =
    \sum_{\s'\in\cell_{j}(M)} [\t:\s'] (\s',\s) = [\t:\s].
    \]
\end{proof}

\begin{prop}
    Let $M$ be a bounded regular CW-complex. Then, for $\sigma\, ,\sigma^{\prime}\in \cell_{j} (M)$, we have
    \[
    (\sigma \, , \Delta_{j+}\, \sigma^{\prime} ) =
    \sum_{\tau\in\cell_{j+1}(M)} [\t:\s] [\t:\s'],
    \]
    and
    \[
    (\sigma \, , \Delta_{j-}\, \sigma^{\prime} ) =
    \sum_{\tau\in\cell_{j-1}(M)} [\s:\t] [\s':\t].
    \]
    In particular,
    \begin{align*}
    (\s, \Delta_{j+} \s) &=
    |\set{\t\in\cell_{j+1}(M) : \dot{\t}\supset\s }|,\\
    (\s, \Delta_{j-} \s) &= |\set{\t\in\cell_{j-1}(M) :
    \t\subset\dot{\s} }|.
    \end{align*}
\end{prop}
\begin{proof}
    Straightforward computation.
\end{proof}

\begin{rem}
    It follows that $\D_{j\pm}$ does not depend on the orientation of
    the $(j\pm 1)$-cells, but only on the orientation of
    the $j$-cells.
\end{rem}

\section{Local Isomorphisms and Geometric Operators}

 In this section, we define geometric operators and prove that the 
 Laplacians (absolute or relative to the boundary subcomplex) are 
 geometric.

\begin{defn}[Combinatorial distance]
    Let $M$ be a connected, regular, bounded
    CW-complex.  Let $\sigma ,\sigma^\prime$ be distinct cells in
    $\cell_j (M)$.  We set

    \itm{i} $d_{-}(\sigma ,\sigma^\prime )=1$, if there is
    $\r\in\cell_{j-1}(M)$ s.t. $\r\subset\dot{\s}\cap \dot{\s}'$,

    \itm{ii} $d_{+}(\sigma ,\sigma^\prime )=1$, if there exists
    $\tau\in \cell_{j+1} (M)$ such that $\s \cup\s' \subset \dot{\t}$,

    \itm{iii} $d(\sigma ,\sigma^\prime )=1$, if either $d_{-}(\sigma
    ,\sigma^\prime )=1$ or $d_{+}(\sigma ,\sigma^\prime )=1$.

    \noindent The distances $d,d_{-},d_{+}$ between two general
    distinct cells $\sigma$ and $\sigma^\prime$ are then defined as
    the minimum number of steps of length one needed to pass from
    $\sigma$ to $\sigma^\prime$, and as $+\infty$ if such a path
    does not exist.
\end{defn}

\noindent We say that $\cell_{j}(M)$ is $d_{\pm}$-connected if
$d_{\pm}(\sigma ,\sigma^\prime ) < +\infty$ for any $\sigma
,\sigma^\prime\in \cell_j (M)$.

\begin{prop} \label{compareDistance}
    Let $M$ be a $p$-dimensional, regular, bounded
    CW-complex.  
    
    \itm{i} If $\cell_{j}(M)$ is $d_{+}$-connected, then it is
    $d_{-}$-connected.  
    
    \itm{ii} Assume any $j$-cell is contained in the boundary of some
    $(j+1)$-cell, $j+1\leq p$.  Then, if $\cell_{j+1}(M)$ is
    $d_{-}$-connected, then $\cell_{j}(M)$ is $d_{+}$-connected.
\end{prop}
\begin{proof}
    $(i)$.  Let us show that if $d_{+}(\s_{0},\s_{1})=1$, $\s_{0}$,
    $\s_{1}\in\cell_{j}(M)$, then $d_{-}(\s_{0},\s_{1})\leq
    V_{j+1}^{-}-1$.  Let $\t\in\cell_{j+1}(M)$ be s.t.
    $\s_{0},\s_{1}\subset\dot{\t}$.  Let $\{\sigma_i\}$ be a basis of
    $j$-cells oriented according to some orientation on $\dot{\t}$,
    which is homeomorphic to the $j$-sphere.  Then $\sum_i \s_i$ is
    the unique $j$-cycle (up to constant multiples) representing the
    non-trivial homology class, hence $\partial_j \sum_i \s_i =0$. 
    This corresponds to the fact that any $(j-1)$-cell has non-trivial
    incidence number with exactly two $j$-cells, one incidence number
    being $1$ and the other $-1$.  Assume now there is a
    $d_{-}$-connected component $\cup_{k}\s_{i_{k}}$ which is properly
    contained in the boundary of $\t$.  Since $\sum_{k}\s_{i_{k}}$ is
    not a cycle, there exists a $(j-1)$-cell $\rho$ such that $(\rho,
    \partial_j \sum_{k} \s_{i_{k}}) \ne 0$.  Then there is exactly one
    $j$-cell, not belonging to $\cup_{k}\s_{i_{k}}$, having
    non-trivial incidence number with $\rho$.  But this is impossible,
    since $\cup_{k}\s_{i_{k}}$ is $d_{-}$-connected.  Since the
    maximum number of $j$-faces of $\t\in\cell_{j+1}(M)$ is $V_{j+1}^{-}$,
    the thesis follows.  \\
    $(ii)$ Let $\r_{1}\ne\r_{2}\in\cell_{j}(M)$,
    $\s_{1},\s_{2}\in\cell_{j+1}$ such that
    $\r_{i}\subset\dot{\s}_{i}$.  Then, since a $d_{-}$-path from
    $\s_{1}$ to $\s_{2}$ gives rise to a $d_{+}$-path from $\r_{1}$ to
    $\r_{2}$, we have
    $$
    d_{+}(\r_{1},\r_{2})\leq d_{-}(\s_{1},\s_{2})+1.
    $$
\end{proof}

 If $\s\in\cell_{j}(M)$, $r\in\bn$, we write $B_{r}(\s) := \set{
 \s'\in\cell_{j}(M): d(\s',\s)\leq r}$.  

\begin{defn}[Finite propagation operators]
    A bounded linear operator $A$ on $\ell^2 (\cell_j M)$ has {\it
    finite propagation} $r=r(A)\geq 0$ if, for all $\sigma\in \cell_j
    (M)$, $\supp (A\sigma)\subset B_{r}(\sigma)$ and $\supp
    (A^{*}\sigma)\subset B_{r}(\sigma)$.
\end{defn}

\begin{lem}
    Finite propagation operators form a $^*$-algebra.
\end{lem}
\begin{proof}
    The set of finite propagation operators is $^*$-closed by
    definition.  To prove that it is also an algebra, one can choose,
    for example,
    $$
    r(\lambda A + B)= r(A)\vee r(B)\, ,\qquad r(AB)= r(A)+ r(B).
    $$
\end{proof}

 Given two CW-complexes $M,\, N$, a continuous map $f:M\to N$ is
 called cellular if $f(M^{j})\subset N^{j}$, for all $j$; it induces
 linear maps $f_{j}: C_{j}(M)\otimes_{\bz}\bc \to
 C_{j}(N)\otimes_{\bz}\bc$ intertwining the boundary maps.  The
 cellular map $f$ is called {\it regular} if, for all $j$,
 $\s\in\cell_{j}(M)$, there are $k$, $\t\in\cell_{k}(N)$ such that
 $f(\s)=\t$, $f(\dot{\s})=\dot{\t}$; then, necessarily, $k\leq j$.  We
 call $f$ an {\it isomorphism} if it is a bijective regular map s.t.
 $[f_{j}\s^{j}_{\a}:f_{j-1}\s^{j-1}_{\b}] =
 [\s^{j}_{\a}:\s^{j-1}_{\b}]$, for all $j,\a,\b$.  Then $f$ is a
 homeomorphism and $f_{j}$ is a linear isomorphism.  
 
 \noindent A subcomplex $N$ of $M$ is a closed subspace of $M$ which
 is a union of (open) cells.  We call $N$ a {\it full} subcomplex if,
 for all $j$, $\s\in\cell_{j}(M)$, $\dot{\s}\subset N$ imply
 $\s\subset N$.

 To prove that a cell belongs to a full subcomplex, we will find it
 convenient in the sequel to refer to the following

 \begin{lem}\label{lem:convenient}
     Let $N$ be a full subcomplex of the regular CW-complex $M$.  Let
     $\t\in\cell_{j}(M)$ be s.t., for all
     $\r\in\cell_{j-1}(\dot{\t})$, one has $\r\in N$.  Then $\t\in N$.
 \end{lem}
 \begin{proof}
     As $N$ is a subcomplex, it follows that $\dot{\t}\subset N$;
     therefore $\t\in N$, because $N$ is full.
 \end{proof}

 \begin{defn}[Local Isomorphisms and Geometric Operators]
     A {\it local isomorphism} of the CW-complex $M$ is a triple
     $$
     \Bigl(s(\gamma)\, ,r(\gamma)\, ,\gamma \Bigr)
     $$
     where $s(\gamma)\, ,r(\gamma)$ are full subcomplexes of $M$ and
     $\gamma : s(\gamma)\to r(\gamma)$ is an isomorphism.
     \vskip0.2truecm
     \noindent For any $j=0,\dots\, ,{\rm dim}(M)$, the
     local isomorphism $\gamma$ defines a {\it partial isometry} $V_j
     (\gamma) : \ell^2 (\cell_j M)\to \ell^2 (\cell_j M)$, by setting
     \begin{align*}
	 V_j(\g)(\s) := 
	 \begin{cases}
	     \g_j(\s) & \s\in\cell_j (s(\g)) \\
	     0 & \s\not\in \cell_j (s(\g)),
	 \end{cases}
     \end{align*}
     and extending by linearity.
     An operator $T\in B(\ell^2 (\cell_j M))$ is called {\it
     geometric} if there exists $r$ such that $T$ has finite
     propagation $r$ and, for any local isomorphism $\gamma$, any
     $\s\in\cell_j(M)$ s.t. $B_{r}(\s)\subset s(\g)$ and
     $B_{r}(\g\s)\subset r(\g)$, one has
     \[
     TV_j (\gamma)\s = V_j (\gamma)T\s,\quad T^*V_j
     (\gamma)\s = V_j (\gamma)T^*\s\, .
     \]
 \end{defn}

 \begin{prop}\label{Prop:3.7}
     Let $M$ be a regular, bounded CW-complex.  Then,
     for any $j$, geometric operators on $\ell^2 (\cell_j M)$ form a
     $^*$-algebra.  The half Laplacians $\Delta_{j\pm}$ belong to it.
 \end{prop}
 \begin{proof}
     The first statement is obvious.  Concerning the second, let us
     note that, since the complex is bounded, half
     Laplacians $\Delta_{j\pm}$ are bounded (cf.  Lemma
     \ref{Lemma:boundedLaplacians}).

     Let $\sigma\, ,\sigma^{\prime}\in \cell_{j} (M)$, with
     $B_{1}(\s)\subset s(\g)$ and $B_{1}(\g\s)\subset r(\g)$.  Then,
     if $\s'\not\in r(\g)$, because $\supp(\D_{j\pm}\s)\subset
     B_{1}(\s) \subset s(\g)$ and $\supp(\D_{j\pm}(\g_{j}\s))\subset
     B_{1}(\g\s) \subset r(\g)$, we get
     $$
     (\s', \D_{j\pm}\, V_{j}(\g)\s) = 0 = (\s',V_{j}(\g)\D_{j\pm}\,
     \s).
     $$
     So, let us suppose that $\s'\in r(\g)$, so that $\s'=\g_{j}\s''$,
     for $\s''\in s(\g)$ and
     \begin{align*}
	 (\s',\D_{j-}\, V_{j}(\g)\s) & = \sum_{\t\in\cell_{j-1}(M)}
	 [\s':\t] [\g_{j}\s:\t] \\
	 & = \sum_{\t'\in\cell_{j-1}(M)} 
	 [\s':\g_{j-1}\t'] [\g_{j}\s:\g_{j-1}\t'] \\
	 & = \sum_{\t'\in\cell_{j-1}(M)} [\s'':\t'] [\s:\t'] \\
	 & = (V_{j}(\g)^{*} \s',\D_{j-}\,\s) = (\s',V_{j}(\g)
	 \Delta_{j-}\, \s),
     \end{align*}
     where the third equality comes from the incidence-preserving
     property of $\g$, and in the second equality we used the fact that
     the non-zero terms in the sum come from $\t$'s which are
     ``components'' of the chain $\de_{j} \g_{j}\s = \g_{j-1}\de_{j}\s
     = \sum c_{i}\g_{j-1}\r_{i}$, if $\de_{j}\s=\sum c_{i}\r_{i}$, so
     that $\t=\g_{j-1}\r_{i}$, for some $i$. By linearity we get that
     $\D_{j-}$ is geometric. As for $\D_{j+}$,
     \begin{align*}
	 (\s',\Delta_{j+}\, V_{j}(\g)\s) & =
	 \sum_{\t\in\cell_{j+1}(M)} [\t:\s'] [\t:\g_{j}\s] \\
	 & = \sum_{\t'\in\cell_{j+1}(M)} [\g_{j+1}\t':\s']
	 [\g_{j+1}\t':\g_{j}\s] \\
	 & = \sum_{\t'\in\cell_{j+1}(M)} [\t':\s''] [\t':\s]  \\
	 & = (V_{j}(\g)^{*} \s',\D_{j+}\,\s) = (\s',V_{j}(\g)
     \Delta_{j+}\, \s),
     \end{align*}
     where the third equality comes from the incidence-preserving
     property of $\g$, and in the second equality we used the fact
     that the non-zero terms in the sum come from $\t$'s s.t.
     $[\t:\g_{j}\s]\neq 0$, so that, for all
     $\r\in\cell_{j}(\dot{\t})$, we get $d(\r,\g\s)=1$, hence $\r\in
     r(\g)$; from Lemma \ref{lem:convenient}, $\t\in r(\g)$, so there
     is $\t'\in s(\g)$ s.t. $\t=\g_{j+1}\t'$.  By linearity we get
     that $\D_{j+}$ is geometric.
 \end{proof}

 We now consider a version of the boundary operators relative to the
 boundary subcomplex.  This idea is due to Lott and L\"{u}ck
 \cite{LoLu}, who introduced relative invariants for covering
 CW-complexes.  In this way, other non-trivial L$^{2}$-Betti numbers are
 available, as shown in Section 8.

 Let $M$ be a $p$-dimensional, regular, bounded CW-complex.  We shall
 consider the $(p-1)$-dimensional {\em boundary subcomplex} $\partial
 M$, defined as follows:
 \begin{itemize}
     \itm{i} a $(p-1)$-cell of $M$ is in $\de M$ if it is
     contained in at most one $p$-cell.
     \itm{ii} a $j$-cell of $M$ is in $\de M$ if it is
     contained in a $(p-1)$-cell in $\de M$.
 \end{itemize}

 Then $\de M$ is a regular bounded CW-complex.

 \begin{lem} \label{lem:scaleDistance} 
     Let $N$ be a full subcomplex of $M$, and $\s_{0}\in\cell_{j}(N)$
     be s.t. $B_{k}(\s_{0})\subset N$.  Then, for any
     $\t_{0}\in\cell_{j+1}(N)$ s.t. $\s_{0}\subset\dot{\t}_{0}$, one
     has $B_{\ell}(\t_{0})\subset N$, for $\ell\leq
     \frac{k}{V_{j+1}^{-}-1}$.
\end{lem}
\begin{proof}
    Let $\t_{1}\in\cell_{j+1}(M)$ be s.t. $d(\t_{1},\t_{0})\leq \ell$. 
    Then, for any $\s_{1}\in\cell_{j}(M)$,
    $\s_{1}\subset\dot{\t}_{1}$, one has, from the proof of
    Proposition \ref{compareDistance}, $d(\s_{1},\s_{0})\leq
    \ell(V_{j+1}^{-}-1)\leq k$.  Therefore $\s_{1}\subset N$.  As $N$
    is full, $\t_{1}\subset N$, and the thesis follows.
\end{proof}

  \begin{lem}\label{Lemma:3.9}
     Let $\g$ be a local isomorphism, $\s\in\cell_{j}(M)$ be s.t.
     $B_{k}(\s)\subset s(\g)$, $B_{k}(\g\s)\subset r(\g)$, where
     $k\geq (V_{p-1}^{-}-1)(V_{p-2}^{-}-1)\ldots(V_{j+1}^{-}-1)$. 
     Then $\s\in \de M$ iff $\g\s\in \de M$.
 \end{lem}
 \begin{proof}
    $(\Rightarrow)$ Let $\s\in\cell_{p-1}(\de M)$.  If there were
    $\t\neq \t'\in\cell_{p}(M)$ s.t.
    $\g\s\subset\dot{\t}\cap\dot{\t'}$, then for all
    $\r\in\cell_{p-1}$, $\r\subset\dot{\t}$ we would get
    $d(\r,\g\s)=1$, hence $\r\in r(\g)$; from Lemma
    \ref{lem:convenient}, $\t\in r(\g)$; analogously $\t'\in r(\g)$. 
    As $\g_{p}$ preserves incidences and boundaries, $\s \subset
    \g^{-1}\dot{\t} \cap \g^{-1}\dot{\t'}$, which implies
    $\s\not\in\de M$, and we have reached a contradiction.  Therefore,
    there is a unique $\t\in\cell_{p}(M)$ s.t. $\g\s \subset
    \dot{\t}$, which means that $\s\in\de M$.

    If $\s\in\cell_{j}(\de M)$, there is $\t\in\cell_{p-1}(M)\cap \de
    M$ s.t. $\s\subset\dot{\t}$, and
    $\g\s\subset(\g\t)\dot{}=\g(\dot{\t})$.  Then, from Lemma
    \ref{lem:scaleDistance}, $B_{1}(\t)\subset\s(\g)$, and
    $B_{1}(\g\t)\subset\r(\g)$.  From what has already been proved,
    $\g\t\in \de M$.  Therefore $\g\s\in \de M$, because $\de M$ is a
    subcomplex.

    $(\Leftarrow)$ follows from the above applied to $\g^{-1}$.
 \end{proof}

 Let $\ov\de_{j} \equiv \de^{M,\de M}_{j}$ be the boundary operator of
 the relative cellular complex $C_{j}(M,\de M) := H_{j}(M^{j}\cup \de
 M, M^{j-1}\cup \de M,\bz)$.  As $\displaystyle C_{j}(M,\de M) \cong
 \oplus_{\s\in\ov\cell_{j}(M)} \bz \s$, where $\ov\cell_{j}(M) :=
 \set{ \s\in\cell_{j}(M) : \s\cap \de M = \emptyset}$, we can identify
 $C^{(2)}_{j}(M,\de M)$, the $\ell^{2}$-completion of $C_{j}(M,\de
 M)\otimes_{\bz}\bc$, with $\ell^{2}(\ov\cell_{j}(M))$, a closed
 subspace of $\ell^{2}(\cell_{j}(M))$.  Moreover we can consider
 $\ov\de_{j}: C^{(2)}_{j}(M)\to C^{(2)}_{j-1}(M)$, $\ov{\de}_{j}^{\,
 *}: C^{(2)}_{j-1}(M)\to C^{(2)}_{j}(M)$, by extending them to $0$ on
 $C^{(2)}_{j}(M,\de M)^{\perp}$ or $C^{(2)}_{j-1}(M,\de M)^{\perp}$,
 respectively.  Define $\ov\D_{j+}:= \ov\de_{j+1}\ov\de_{j+1}^{\, *}$,
 $\ov\D_{j-}:= \ov\de_{j}^{\, *}\ov\de_{j}$.  Then

 \begin{lem}
      \itm{i} $\ov\D_{j\pm}\s = 0$, for $\s\in C^{(2)}_{j}(M,\de
      M)^{\perp}$, 
      
      \itm{ii} for $\s,\s'\in C^{(2)}_{j}(M,\de M)$,
      \[
      (\s'\, , \ov\D_{j+}\, \s) = \sum_{\t\in\ov\cell_{j+1}(M)}
      [\t:\s] [\t:\s']\, , \qquad (\s'\, , \ov\D_{j-}\, \s) =
      \sum_{\t\in\ov\cell_{j-1}(M)} [\s:\t] [\s':\t]\, .
      \]
 \end{lem}


 \begin{prop} 
     Let $\g$ be a local isomorphism, $\s\in\cell_{j}(M)$ be such that
     $B_{k}(\s)\subset s(\g)$, $B_{k}(\g\s)\subset r(\g)$, for some
     $k\geq 1+ \prod_{i=j+1}^{p-1}(V_{i}^{-}-1)$.  Then
     $$
     \ov\D_{j-}\, V_{j}(\g) \s = V_{j}(\g)\, \ov\D_{j-} \s, \qquad
     \ov\D_{j+}\, V_{j}(\g) \s = V_{j}(\g)\, \ov\D_{j+} \s\, .
     $$
 \end{prop}
 \begin{proof}
     Let us prove that, for any $\s'\in\ov\cell_{j}(M)$, we have
     $$
     (s',\ov\D_{j\pm}\, V_{j}(\g) \s) = (\s',V_{j}(\g)\,
     \ov\D_{j\pm}\s).
     $$
     If $\s'\not\in B_{1}(\g\s)$, the thesis is true.  Indeed, from
     $\supp(\ov\D_{j\pm}\s)\subset B_{1}(\s)\subset s(\g)$, and
     $\supp(\ov\D_{j\pm}\g\s)\subset B_{1}(\g\s)\subset r(\g)$, it
     follows $(\s',\ov\D_{j\pm} V_{j}(\g)\s) = 0$, whereas, if
     $\s'\not\in r(\g)$ we get $(\s',V_{j}(\g)\ov\D_{j\pm} \s) = 0$,
     while, if $\s'\in r(\g)\setminus B_{1}(\g\s)$, we get
     $(\s',V_{j}(\g)\ov\D_{j\pm} \s) = (\g_{j}^{-1}\s',\ov\D_{j\pm}
     \s) = 0$, as $d(\g_{j}^{-1}\s',\s) = d(\s',\g_{j}\s) >1$. 
     Therefore, we can assume $\s'\in B_{1}(\g\s)$.  Moreover, if
     $\s\in\de M$, so that $\g\s\in\de M$ (by Lemma \ref{Lemma:3.9}),
     we get $\ov\D_{j\pm} V_{j}(\g)\s = 0 = V_{j}(\g)\ov\D_{j\pm} \s$.

     Therefore, we now assume $\s\not\in\de M$, $\s'\in B_{1}(\g\s)$.
     Then 
     $$
     (\s',\ov\D_{j-}\, V_{j}(\g) \s) = \sum_{\t\in\ov\cell_{j-1}(M)}
     [\s':\t][\g_{j}\s:\t].
     $$
     Let
     $\t\in\ov\cell_{j-1}(M)$, $\t\subset (\g\s)\dot{}\cap\dot{\s}'$.
     Then, as in the first part of the proof of Proposition
     \ref{Prop:3.7}, there is $\t'\in\cell_{j-1}(s(\g))$ s.t.
     $\t=\g\t'$; moreover $\t'\subset\dot{\s}$, as $[\s:\t'] =
     [\g\s:\g\t'] \neq 0$.  Let us now show that $\t\in\de M \iff
     \t'\in\de M$; indeed, if $\t\in\de M$, then there is
     $\r\in\cell_{j}(\de M)$ s.t. $\t\subset\dot{\r}$; therefore
     $d(\r,\g\s)\leq1$, and $B_{k-1}(\r)\subset r(\g)$, so, from Lemma
     \ref{Lemma:3.9}, it follows that $\r':= \g_{j}^{-1}\r\in\de M$;
     then $[\r':\t'] = [\r:\t] \neq 0$, hence $\t'\subset\dot{\r}'$,
     and $\t'\in\de M$.  The other implication follows similarly.
     Therefore
     \begin{align*}
	 (\s',\ov\D_{j-}\, V_{j}(\g) \s) & =
	 \sum_{\t\in\ov\cell_{j-1}(M)} [\s':\t][\g_{j}\s:\t] \\
	 & = \sum_{\t'\in\ov\cell_{j-1}(M)}
	 [\s':\g_{j-1}\t'] [\g_{j}\s:\g_{j-1}\t'] \\
	 & = \sum_{\t'\in\ov\cell_{j-1}(M)}
	 [\g_{j}^{-1}\s':\t'][\s:\t'] \\
	 & = (V_{j}(\g)^{*} \s',\ov\D_{j-}\s) = (\s',V_{j}(\g)\,
	 \ov\D_{j-}\s).
     \end{align*}

     As for $\ov\D_{j+}$, we get 
     $$
     (\s',\ov\D_{j+}\, V_{j}(\g) \s)
     = \sum_{\t\in\ov\cell_{j+1}(M)} [\t:\s'][\t:\g_{j}\s].
     $$
     Let $\t\in\cell_{j+1}(M)$ be s.t. $\g\s\cup\s' \subset\dot{\t}$;
     as for any $\r\in\cell_{j}(M)$, $\r\subset\dot{\t}$, it holds
     $d(\r,\g\s)\leq 1$, so $\r\in B_{1}(\g\s)\subset r(\g)$, from
     Lemma \ref{lem:convenient} we get $\t\in r(\g)$; therefore there
     is $\t'\in s(\g)$ s.t. $\t=\g\t'$.  From Lemma
     \ref{lem:scaleDistance} it follows $B_{\ell}(\t)\subset r(\g)$,
     for $\ell\leq k/(V_{j+1}^{-}-1)$, and Lemma \ref{Lemma:3.9} gives
     us $\t\in\de M \iff \t'\in\de M$.  Therefore
     \begin{align*}
	 (\s',\ov\D_{j+}\, V_{j}(\g) \s) & =
	 \sum_{\t\in\ov\cell_{j+1}(M)} [\t:\s'][\t:\g_{j}\s] \\
	 & = \sum_{\t'\in\ov\cell_{j+1}(M)}
	 [\g_{j+1}\t':\s'][\g_{j+1}\t':\g_{j}\s] \\
	 & = \sum_{\t'\in\ov\cell_{j+1}(M)}
	 [\t':\g_{j}^{-1}\s'][\t':\s] \\
	 & = (V_{j}(\g)^{*} \s',\ov\D_{j+}\s) = (\s',V_{j}(\g)\,
	 \ov\D_{j+}\s).
     \end{align*}
 \end{proof}

 We have proved the following.

 \begin{prop}
     Let $M$ be a $p$-dimensional, regular, bounded
     CW-complex.  The relative half-Laplacians $\overline{\D}_{j\pm}$
     are geometric operators.
 \end{prop}

\section{Self-similar CW-complexes}

 In this section we introduce self-similar complexes, and show that
 there is a natural trace state on the algebra of geometric operators.

 If $K$ is a subcomplex of $M$, we call $j$-{\it frontier} of $K$, and
 denote it by $\cf(\cell_{j}K)$, the family of cells in $\cell_{j}K$
 having distance 1 from the complement of $\cell_{j}K$ in $\cell_{j}(M)$.
 
 \begin{defn}[Amenable CW-Complexes]
     A countably infinite CW-complex $M$ is {\it amenable} if it is
     regular and bounded, and has an {\it amenable
     exhaustion}, namely, an increasing family of finite subcomplexes
     $\{K_n : n\in \mathbb{N}\}$ such that $\cup K_n = M$ and for all
     $j=0,\dots\, ,{\rm dim}(M)$,
    \begin{equation*}
	\frac{|\cf(\cell_j K_n)|}{|\cell_j K_n|}\to 0\qquad {\rm
	as}\,\,\, n\to \infty\, .
    \end{equation*}
 \end{defn}

 \begin{defn}[Self-similar CW-Complexes] \label{def:Quasiperiodic} 
     A countably infinite CW-complex $M$ is {\it self-similar} if it
     is regular and bounded, and it has an amenable exhaustion by full
     subcomplexes $\set{K_{n} : n\in\bn}$ such that the following 
     conditions $(i)$ and $(ii)$ hold:
     
     \itm{i} for all $n$ there is a finite set of local isomorphisms
     $\cg(n,n+1)$ such that, for all $\gamma\in \cg(n,n+1)$, one has
     $s(\gamma) = K_n$,
     \begin{equation*}
	 \bigcup_{\gamma\in \cg(n,n+1)} \gamma_j \Bigl(\cell_j (K_n)
	 \Bigr)= \cell_j (K_{n+1} ),\qquad j=0,\dots\, ,{\rm dim}(M)
     \end{equation*}
     and moreover if $\gamma , \gamma^\prime\in \cg(n,n+1)$ with
     $\gamma\neq\gamma^\prime$
     \begin{equation}\label{1BdryIntersec}
	 \cell_j \gamma(K_n ) \cap \cell_j \gamma'(K_n ) = \cf(\cell_j
	 \gamma(K_n )) \cap \cf(\cell_j \gamma'(K_n )) \, ,\quad
	 j=0,\dots\, ,{\rm dim}(M).
     \end{equation}
    
     \itm{ii} We then define $\cg(n,m)$, with $n<m$, as the set of all
     admissible products $\g_{m-1}\cdot\dots\cdot\g_{n}$, $\g_{i}\in
     \cg(i,i+1)$, where admissible means that the range of $\g_{j}$ is
     contained in the source of $\g_{j+1}$.  We let $\cg(n,n)$ consist
     of the identity isomorphism on $K_{n}$, and $\cg(n)=\cup_{m\geq
     n}\cg(n,m)$.  We now define the $\cg$-invariant $j$-frontier of
     $K_{n}$:
     $$
     \cf_{\cg}(\cell_{j}K_{n})= \bigcup_{\gamma\in
     \cg(n)}\g_{j}^{-1}\cf(\cell_{j}\gamma(K_{n})),
     $$
     and we ask that
     \begin{equation*}
	 \frac{|\cf_{\cg}(\cell_j K_n)|}{|\cell_j K_n|}\to 0\qquad
	 {\rm as}\; n\to \infty\, .
     \end{equation*}
 \end{defn}
 
 \begin{rem}\label{rem:rBdryIntersec}
     We may replace the condition in (\ref{1BdryIntersec}) with
     the following
     $$
     \cell_j \gamma(K_n ) \cap \cell_j \gamma'(K_n ) \subseteq B_{r}
     (\cf(\cell_j \gamma(K_n ))) \cap B_{r}(\cf(\cell_j \gamma'(K_n
     ))) \, ,\quad j=0,\dots\, ,{\rm dim}(M),
     $$
     for a suitable $r>0$.  It is easy to see that all the theory
     developed below will remain valid.
 \end{rem}
 
 Some examples of self-similar CW-complexes are given below, cf.
 Section \ref{sec:prefractals} for more details on the construction.

 \begin{ex}
    The Gasket graph in figure \ref{fig:Gasket}, the Lindstrom
    graph in figure \ref{fig:Lindstrom}, the Vicsek graph in figure
    \ref{fig:Vicsek} are examples of $1$-dimensional self-similar 
    complexes. The Carpet $2$-complex in figure \ref{fig:Carpet}
    is an example of a $2$-dimensional self-similar CW-complex.
    \begin{figure}[ht]
    \centering
    \psfig{file=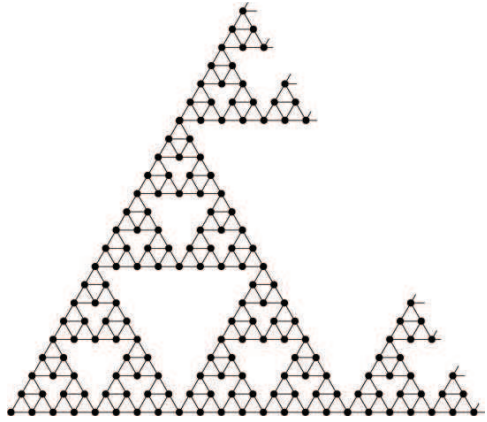,height=2.165in,width=2.5in}
    \caption{Gasket graph}
    \label{fig:Gasket}
    \end{figure}

    \begin{figure}[ht]
    \centering
    \psfig{file=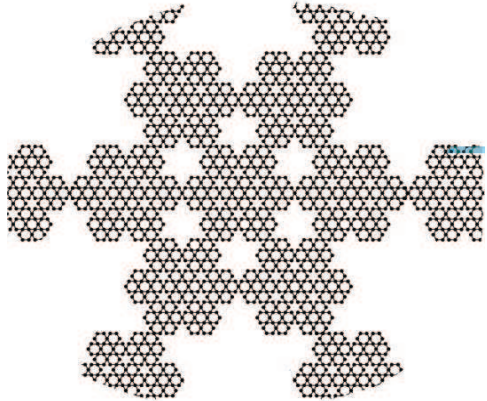,height=2.165in,width=2.5in}
    \caption{Lindstrom graph}
    \label{fig:Lindstrom}
    \end{figure}

    \begin{figure}[ht]
    \centering
    \psfig{file=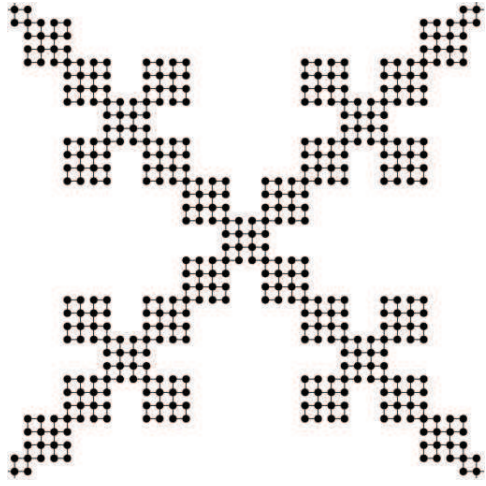,height=2.5in,width=2.5in}
    \caption{Vicsek graph}
    \label{fig:Vicsek}
    \end{figure}

    \begin{figure}[ht]
    \centering
    \psfig{file=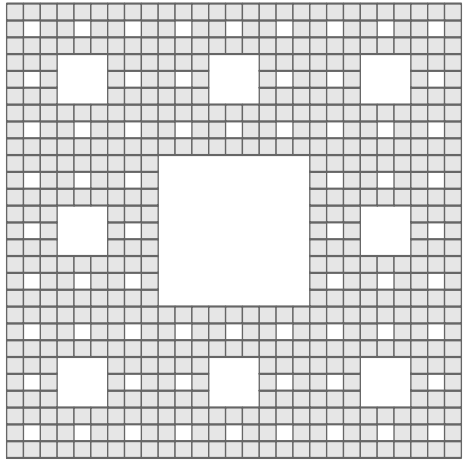,height=2.5in,width=2.5in}
    \caption{Carpet $2$-complex}
    \label{fig:Carpet}
    \end{figure}
 \end{ex}

 \begin{thm}\label{TraceThm}
    Let $M$ be a self-similar CW-complex, $\ca (\cell_j M)$ the
    C$^*$--algebra given by the closure of the $^{*}$-algebra of
    geometric operators.  Then, on $\ca(\cell_{j}M)$ there is a well defined
    trace state $\Phi_{j}$ given by
    \begin{equation*}
	\Phi_{j} (T) = \lim_n \frac{Tr\bigl( E(\cell_j K_n) T \bigr) }
	{Tr\bigl( E(\cell_j K_n) \bigr) }
    \end{equation*}
    where $E(\cell_j K_n)$ is the orthogonal projection of
    $\ell^2(\cell_j M)$ onto $\ell^{2}(\cell_j K_n)$.
 \end{thm}
 \begin{proof}
    Fix $j\in\set{0,\ldots,p}$, and for a finite subset $N\subset
    \cell_{j}M$ denote by $E(N)\in \cb(\ell^2 (\cell_j M))$ the
    projection onto $\text{span} N$.  Let us observe that, since $N$
    is an orthonormal basis for $\ell^2 (N)$, then $Tr \bigl(E(N)
    \bigr) = |N|$.  
    
    \noindent {\bf First step}: some combinatorial results.
    
    a) Let $\m\equiv\m_{j}=\sup_{\s\in\cell_{j}M}|B_{1}(\s)|$.  First
    observe that $\m$ is finite, since $\m\leq V^{+}_{j}+V^{-}_{j}$.
    
    Then, since
    $$B_{r+1}(\s)=\bigcup_{\s'\in B_{r}(\s)}B_{1}(\s'),$$
    we get $|B_{r+1}(\s)|\leq|B_{r}(\s)|\m$, giving
    $|B_{r}(\s)|\leq \m^{r}$, $\forall\s\in\cell_{j}M$, $r\geq 0$.
    As a consequence, for any finite set $\Omega\subset\cell_{j}M$,
    we have
    $B_{r}(\Omega)=\cup_{\s'\in \Omega}B_{r}(\s')$,
    giving
    \begin{equation}\label{subsets}
    |B_{r}(\Omega)|\leq |\Omega|\m^{r},\quad\forall r\geq 0.
    \end{equation}

    b) Let us set $\O(n,r)=\cell_{j}K_{n}\setminus
    B_{r}(\cf_{\cg}(\cell_{j}K_{n}))$.  Then, for any $\g\in\cg(n)$,
    we have
    \begin{equation*}
	\g_{j}\O(n,r)\subset\g_{j}\cell_{j}K_{n}\subset\g_{j}\O(n,r)
	\cup B_{r}(\cf_{\cg}(\g_{j}\cell_{j} K_{n})).
    \end{equation*}
    Now assume $r\geq1$.  Then, the $\g_{j}\O(n,r)$'s are disjoint, for
    different $\g$'s in $\cg(n,m)$.  Therefore,
    \begin{align}
	|\cell_{j}K_{n}| &\leq |\O(n,r)| + |\cf_{\cg}(\cell_{j}K_{n})|
	\m^{r}, \label{ineqA}\\
	\left|\cell_{j}K_{m}\setminus\bigcup_{\g\in\cg(n,m)}
	\g_{j}\O(n,r)\right| &\leq
	|\cg(n,m)|\,|\cf_{\cg}(\cell_{j}K_{n})| \m^{r},\label{ineqB}\\
	|\cg(n,m)|\,|\O(n,r)| & \leq |\cell_{j}K_{m}| \leq
	|\cg(n,m)|\,|\cell_{j}K_{n}|.\label{ineqC}
    \end{align}
    Indeed, (\ref{ineqA}) and (\ref{ineqC}) are easily verified, whereas
    \begin{align*}
	\left|\cell_{j}K_{m}\setminus\bigcup_{\g\in\cg(n,m)}
	\g_{j}\O(n,r)\right| & = \left|\bigcup_{\g\in\cg(n,m)}\g_{j}
	\cell_{j}K_{n}\setminus\bigcup_{\g\in\cg(n,m)}\g_{j}\O(n,r)\right|
	\\
	& \leq \sum_{\g\in\cg(n,m)} \left|\g_{j}
	[\cell_{j}K_{n}\setminus \O(n,r)]\right| \\
	& \leq |\cg(n,m)| \left|
	B_{r}(\cf_{\cg}(\cell_{j}K_{n}))\right| \\
	& \leq |\cg(n,m)| \left| \cf_{\cg}(\cell_{j}K_{n})\right|
	\m^{r}.
    \end{align*}

    c) Let $\eps_{n}=\displaystyle{\frac{|\cf_{\cg}(\cell_{j}K_{n})|}
    {|\cell_{j}K_{n}|}}$, and recall that $\eps_{n}\to0$.  Putting
    together (\ref{ineqA}) and (\ref{ineqC}) we get
    $$
    |\cg(n,m)|\,|\cell_{j}K_{n}| -
    |\cg(n,m)|\,|\cf_{\cg}(\cell_{j}K_{n})| \m^{r} \leq
    |\cell_{j}K_{m}| \leq |\cg(n,m)|\,|\cell_{j}K_{n}|,
    $$
    which implies
    $$
    1-\eps_{n}\m^{r}\leq \frac{|\cell_{j}K_{m}|} {|\cg(n,m)|\,
    |\cell_{j}K_{n}|} \leq 1.
    $$
    Choosing $n_{0}$ such that, for $n>n_{0}$,
    $\eps_{n}\m^{r}\leq1/2$, we obtain
    \begin{equation}\label{ineqE}
	0\leq \frac{|\cg(n,m)|\, |\cell_{j}K_{n}|}{|\cell_{j}K_{m}|}
	-1 \leq 2\eps_{n}\m^{r}\leq 1.
    \end{equation}
    Therefore, from (\ref{ineqB}), we obtain
    \begin{align}\label{ineqF}
	\left|\cell_{j}K_{m}\setminus\bigcup_{\g\in\cg(n,m)}
	\g_{j}\O(n,r)\right| & \leq |\cg(n,m)| \left|
	\cf_{\cg}(\cell_{j}K_{n})\right| \m^{r} \\
	& = |\cg(n,m)| \left| \cell_{j}K_{n}\right| \eps_{n}\m^{r}
	\leq 2 \left| \cell_{j}K_{m}\right| \eps_{n}\m^{r}.  \notag
    \end{align}
    \\
    {\bf Second step}: the existence of the limit for geometric
    operators.
    
    a) By definition of $V_j (\gamma )$, we have,
    for $\gamma\in \cg(n,m)$, $n<m$,
    $$
    V_j^* (\gamma) V_j (\gamma) = E(\cell_j K_n), \qquad V_j (\gamma)
    V_j^* (\gamma) = E\Bigl(\gamma_j \bigl(\cell_j K_n \bigr)\Bigr).
    $$
    Assume now $T\in \mathcal{B}(\ell^2 (\cell_j M))$ is a geometric
    operator with finite propagation $r$.  Then,
    $$
    T V_j(\gamma) E(\O(n,r))= V_j(\gamma) T E(\O(n,r)) \qquad
    E(\gamma_j \O(n,r)) = V_j (\gamma) E(\O(n,r)) V_j (\gamma)^*.
    $$
    As a consequence,
    \begin{equation}\label{eq:5.6}
    \begin{split}
	Tr\bigl(T& E(\gamma_j(\O(n,r))) \bigr) \\
	&= Tr\bigl(T V_j(\gamma) E(\O(n,r)) V_j (\gamma)^* \bigr) =
	Tr\bigl(V_j(\gamma) T E(\O(n,r)) V_j (\gamma)^* \bigr)\\
	&= Tr \bigl(TE(\O(n,r)) V_j (\gamma)^* V_j (\gamma) \bigr) =
	Tr \bigl(TE(\O(n,r)) E(\cell_j K_n) \bigr)\\
        &=Tr \bigl(TE(\O(n,r)) \bigr) .
    \end{split}
    \end{equation}
    
    b) Let us show that the sequence is Cauchy:
    \begin{align*}
	&\left|\frac{Tr TE(\cell_{j}K_{n})}{Tr E(\cell_{j}K_{n})} -
	\frac{Tr TE(\cell_{j}K_{m})}{Tr E(\cell_{j}K_{m})} \right| \\
	\leq & \frac{|Tr T(E(\cell_{j}K_{n}) - E(\O(n,r)))|}
	{|\cell_{j}K_{n}|} + \frac {|Tr T(E(\cell_{j}K_{m}) -
	E(\cup_{\g\in\cg(n,m)}\g_{j}\O(n,r)))|} {|\cell_{j}K_{m}|} \\
	+ & \left|\frac{Tr T E(\O(n,r))}{|\cell_{j}K_{n}|} -
	\frac{|\cg(n,m)|\, |\cell_{j}K_{n}|}{|\cell_{j}K_{m}|}
	\frac{Tr T E(\O(n,r))}{|\cell_{j}K_{n}|} \right|\\
	\leq & \|T\|\left( \frac{|\cell_{j}K_{n} \setminus
	\O(n,r)|}{|\cell_{j}K_{n}|} + \frac{|\cell_{j}K_{m} \setminus
	\cup_{\g\in\cg(n,m)}\g_{j}\O(n,r)|} {|\cell_{j}K_{m}|} +
	\left|1 - \frac{|\cg(n,m)|\,
	|\cell_{j}K_{n}|}{|\cell_{j}K_{m}|} \right|\right)\\
	\leq & 5\|T\|\eps_{n}\m^{r},
    \end{align*}
    where we used (\ref{eq:5.6}), in the first inequality, and
    (\ref{ineqF}), (\ref{ineqE}), in the second inequality.  
    \medskip
    
    \noindent {\bf Third step}: $\Phi_{j}$ is a state on $\ca(\cell_j(M))$.
    
    a) Let $T\in\ca(\cell_j(M))$, $\eps>0$. Now find a geometric
    operator $T'$ such that $\|T-T'\|\leq\eps/3$, and set $\phi_{n}(A)
    := \frac{Tr AE(\cell_{j}K_{n})} {Tr E(\cell_{j}K_{n})}$.  Then
    choose $n$ such that, for $m>n$, $|\phi_{m}(T') - \phi_{n}(T')|
    \leq \eps/3$. We get
    $$
    |\phi_{m}(T)-\phi_{n}(T)|\leq |\phi_{m}(T-T')| +
    |\phi_{m}(T')-\phi_{n}(T')| + |\phi_{n}(T-T')|\leq\eps
    $$
    namely $\lim \phi_{n}(T)$ exists.
    
    b) The functional $\Phi_{j}$ is clearly linear, positive and
    takes value $1$ at the identity, hence it is a state on $\ca
    (\cell_j (M))$.
    \medskip
    
    \noindent {\bf Fourth step}: $\Phi_{j}$ is a trace on $\ca(\cell_j(M))$.
    
    Let $A$ be a geometric operator with propagation $r$.  Then
    \begin{align*}
	AE(\cell_{j}K_{n}) &=
	E(B_{r}(\cell_{j}K_{n}))AE(\cell_{j}K_{n}), \\
	E(\O(n,r))A &= E(\O(n,r)) A E(\cell_{j}K_{n}).
    \end{align*}
    Indeed, 
    $$
    \O(n,r) \subset \cell_{j}K_{n} \setminus
    B_{r}(\cf(\cell_{j}K_{n})) = \set{\s\in\cell_{j}K_{n}:
    d(\s,M\setminus\cell_{j}K_{n})\geq r+2},
    $$
    so that 
    $$
    B_{r}(\O(n,r)) \subset \set{\s\in\cell_{j}K_{n}:
    d(\s,M\setminus\cell_{j}K_{n})\geq 2} \subset \cell_{j}K_{n}.
    $$
    Since $A^{*}$ has propagation $r$, we get 
    $$
    A^{*}E(\O(n,r)) = E(B_{r}(\O(n,r))A^{*}E(\O(n,r)) =
    E(\cell_{j}K_{n})A^{*}E(\O(n,r)),
    $$
    which proves the claim.  Therefore,
    \begin{align*}
	&AE(\cell_{j}K_{n}) = E(B_{r}(\cell_{j}K_{n}) \setminus
	\O(n,r))AE(\cell_{j}K_{n}) + E(\O(n,r))A \\
	&= E(B_{r}(\cell_{j}K_{n}) \setminus
	\O(n,r))AE(\cell_{j}K_{n}) - E(\cell_{j}K_{n} \setminus
	\O(n,r))A + E(\cell_{j}K_{n}) A.
    \end{align*}
    Therefore, if $B\in\ca (\cell_j (M))$,
    \begin{align*}
	\phi_{n}([B,A])&\leq\|A\|\,\|B\| \frac{|B_{r}(\cell_{j}K_{n})
	\setminus \O(n,r)|+ |\cell_{j}K_{n} \setminus
	\O(n,r)|}{|\cell_{j}K_{n}|}\\
	&\leq 2 \|A\|\,\|B\|\eps_{n}\m^{r},
    \end{align*}
    as $B_{r}(\cell_{j}K_{n}) \setminus \O(n,r) \subset
    B_{r}(\cf_{\cg}(\cell_{j}K_{n}))$.  Taking the limit as
    $n\to\infty$ we get $\Phi_{j}([B,A])=0$.  By continuity, the
    result holds for any $A,B\in\ca (\cell_j (M))$.
\end{proof}

In the following we use a different normalisation for the
traces and, by giving up the state property, we obtain that  the trace of the identity operator 
in $\ca_{j}$ measures the relative volume of $\cell_{j}(M)$.
This simplifies the relations in Corollaries \ref{NSequality} and  
\ref{EulerPoincare}.

\begin{lem}\label{lem:diffNorm}
    Let $M$ be a $p$-dimensional self-similar complex.  The following
    limits exist and are finite:
    $$
    \lim_{n}\frac{|\cell_{j}(K_{n})|}{|\cell_{p}(K_{n})|}, \quad 0\leq
    j\leq p.
    $$
\end{lem}
\begin{proof}
    We show that the sequences are Cauchy.  Indeed, by inequalities
    (\ref{ineqE}) in the proof of Theorem \ref{TraceThm}, we have, for
    $m>n$ and $j=0,\dots,p$,
    $$
    (1+2\eps_{n}\mu)^{-1}|\cg(n,m)| |\cell_{j}K_{n}|
    \leq |\cell_{j}K_{m}|
    \leq |\cg(n,m)| |\cell_{j}K_{n}|,
    $$
    where the sequence $\eps_{n}=\sup_{j=1,\dots,p}
    \frac{|\cf_{\cg}(\cell_{j}K_{n})|} {|\cell_{j}K_{n}|}$ 
    is infinitesimal and less than 1, and 
    $\mu=\sup_{j}\sup_{\s\in\cell_{j}M}|B_{1}(\s)|$. Therefore 
    \begin{equation*}\label{ineqZ}
	(1-2\eps_{n}\mu)\frac{|\cell_{j}K_{n}|}{|\cell_{p}K_{n}|}
	\leq \frac{|\cell_{j}K_{m}|}{|\cell_{p}K_{m}|}
	\leq (1+2\eps_{n}\mu)\frac{|\cell_{j}K_{n}|}{|\cell_{p}K_{n}|}.
    \end{equation*}
    Hence, the sequence 
    $\frac{|\cell_{j}K_{n}|}{|\cell_{p}K_{n}|}$ is bounded by 
    some constant $M>0$, and   
    $$
    \left|\frac{|\cell_{j}K_{m}|}{|\cell_{p}K_{m}|} -
    \frac{|\cell_{j}K_{n}|}{|\cell_{p}K_{n}|}\right| \leq
    2M\mu\eps_{n}.
    $$
    The thesis follows.
\end{proof}

\begin{defn}\label{Def:BNS}
    Let $M$ be a $p$-dimensional self-similar complex. On the 
    C$^{*}$-algebras $\ca_{j}$ we shall consider the traces
    $$
    Tr^{\cg}_{j}(T)=\lim_{n}\frac{|\cell_{j}(K_{n})|}
    {|\cell_{p}(K_{n})|} \Phi_{j}(T) = \lim_n \frac{Tr\bigl( E(\cell_j
    K_n) T \bigr) } {Tr\bigl( E(\cell_p K_n) \bigr) }.
    $$
    In this way, $Tr^{\cg}_{j}(I)$ measures the relative volume of
    $\cell_{j}(M)$ w.r.t.  $\cell_{p}(M)$.
\end{defn}

\section{L$^2$-Betti numbers and Novikov-Shubin numbers for
self-similar CW-complexes}

In this section, we define L$^{2}$-Betti numbers and Novikov-Shubin numbers for
self-similar CW-complexes, prove various relations among them, and 
give a result on the Euler-Poincar\'e caracteristic of a complex.

Let $M$ be a self-similar CW-complex, let $\D$ be one of the 
operators $\D_{j\pm}$, $\D_{j}$, $\ov{\D}_{j\pm}$, $\ov{\D}_{j}$, and 
define
\begin{defn}[L$^2$-Betti and Novikov-Shubin numbers]
    \itm{i} $\displaystyle{ \b(\D) := \lim_{t\to\infty}
    Tr^{\cg}_{j}(e^{-t\D}), }$ the L$^{2}$-Betti number of $\D$,
    
    \itm{ii} $\displaystyle{ \a(\D) := 2
    \lim_{t\to\infty} \frac{\log \left( Tr^{\cg}_{j}(e^{-t\D}) -
    \b(\D) \right)} {-\log t} }$, the Novikov-Shubin number of $\D$,
    
    and the lower and upper versions, if the above limits do not 
exist.
\end{defn}

Then set
    \begin{align*}
    \b_{j}^{\pm}(M) & := \b(\D_{j\pm}),  \\
    \b_{j}(M) & := \b(\D_{j}), \\
    \b_{j}^{\pm}(M,\partial M) & := \b(\overline{\D}_{j\pm}),\\
    \b_{j}(M,\partial M) & := \b(\overline{\D}_{j}),
    \end{align*}
and analogously for the Novikov-Shubin numbers.

\begin{rem}
    $(i)$ The $L^{2}$-Betti numbers and Novikov-Shubin numbers could
    have been defined also in terms of the spectral density function
    $N_{\lambda}$.  This is usually defined in terms of spectral
    projections, which belong to the generated von Neumann algebra,
    hence, in our case, are not necessarily in the domain of the
    trace.  However we may consider the spectral measure $\m_{j\pm}$
    associated, via Riesz theorem, to the functional $\varphi\in
    C^{0}[0,\infty)\mapsto Tr^{\cg}_{j} (\varphi(\D_{j\pm}))\in \bc$,
    and then define $N_{\lambda}(\D_{j\pm}) :=\int_{0}^{\l}
    d\m_{j\pm}$.  The two definitions for the L$^2$-Betti numbers
    clearly coincide, while Novikov-Shubin numbers can be related via
    a Tauberian theorem, as in \cite{GS}.
    
    $(ii)$ We followed \cite{LoLu} for the definition of the
    relative $L^{2}$-invariants, even though we considered only the
    two cases of no boundary and of full boundary.  It would be
    interesting to prove their Poincar\'e duality result in our
    context.
    
    $(iii)$ We have used the same normalizing sequence for each 
    trace $Tr_j^{\cg}$, see Definition \ref{Def:BNS}, in order to 
    compare $L^{2}$-Betti numbers. This will imply the relation in 
    Corollary \ref{EulerPoincare}.
\end{rem}

\begin{lem}[Hodge decomposition]
    The following decomposition holds true:
    \[
    \ell^2 \bigl(\cell_j M\bigr) = \overline{{\rm Im} \Delta_{j-}}
    \oplus \overline{{\rm Im} \Delta_{j+}} \oplus {\rm ker} 
\Delta_j\, .
    \]
\end{lem}
\begin{proof}
    \[
    \overline{{\rm Im} \Delta_{j+}} = \bigl({\rm ker}\,
    \partial_{j+1}\, \partial_{j+1}^* \bigr)^\perp = \bigl({\rm
    ker}\, \partial_{j+1}^* \bigr)^\perp = \overline{{\rm Im}\,
    \partial_{j+1}} \subseteq {\rm ker}\, \partial_j = {\rm
    ker}\,\partial_j^*\, \partial_j = \bigl( {\rm Im}
    \Delta_{j-}\bigr)^\perp\, .
    \]
    Then since
    \[
    \bigl( {\rm Im} \Delta_{j+}\bigr)^\perp \cap \bigl( {\rm Im}
    \Delta_{j-}\bigr)^\perp = {\rm ker}\, \partial_j \cap {\rm
    ker}\, \partial_{j+1}^*  = {\rm ker}\, \Delta_j
    \]
    the thesis follows.
\end{proof}

\begin{thm}\label{thm:BettiAndNS}
    With the notation above, we have the relations:
    \begin{align*}
	\beta_j(M) &= \beta_{j}^{+}(M) + \beta_{j}^{-}(M)
	-Tr^{\cg}_{j}(I),\\
	\beta_j(M,\partial M) &= \beta_{j+}(M,\partial M) +
	\beta_{j-}(M,\partial M) -Tr^{\cg}_{j}(I),\\
	\alpha_j(M) &= \min\{\alpha_{j+}(M)\, , \alpha_{j-}(M)\},\\
	\alpha_j(M,\partial M) &= \min \{\alpha_{j+}(M,\partial M)\,
	, \alpha_{j-}(M,\partial M)\}.
    \end{align*}
\end{thm}
\begin{proof}
    By the orthogonality of the ranges we have $\Delta_{j+}\,
    \Delta_{j-} = 0$.  Hence $\bigl(\Delta_{j+} + \Delta_{j-}\bigr)^n
    = \Delta_{j+}^n + \Delta_{j-}^n$ from which we get $e^{-t\Delta_j}
    = e^{-t\Delta_{j+}} + e^{-t\Delta_{j-}} - I$.  Now
    the thesis easily follows.  The proof for the relative invariants
    is analogous.
\end{proof}

\begin{prop}
    $\displaystyle{ \Tr_{j-1}((\de_{j}\de_{j}^{*})^{k}) = 
    \Tr_{j}((\de_{j}^{*}\de_{j})^{k}) , \ k\in\bn}$
\end{prop}
\begin{proof}
    Let us set $\O_n := \cell_{j-1} K_n \setminus
    \cf_\cg(\cell_{j-1}K_n)$.  First we note that 
    $$
    E(\O_n) \de_{j} = E(\O_n) \de_{j} E(\cell_{j}K_{n}),$$ indeed
    they coincide on the range of $E(\cell_{j}K_{n})$, and both vanish
    on its kernel.  Analogously $\displaystyle{ \de_{j}
    E(\cell_{j}K_{n}) = E(B_1(\cell_{j-1}K_{n})) \de_{j}
    E(\cell_{j}K_{n}) }$.  \\
    Let us note that if $\de_{j} = V_{j}
    |\de_{j}|$ denotes the polar decomposition, and
    $A\in\cb(\ell^{2}(\cell_{j}(M)))$, then
    $\displaystyle{ Tr( \de_{j}A\de_{j}^{*}) = Tr(
    |\de_{j}|A|\de_{j}|) }$. Then
    $$
    Tr(E(\O_n) (\de_{j}\de_{j}^{*})^{k}) = Tr( E(\O_n) \de_{j}
    E(\cell_{j}K_{n}) \de_{j}^{*} (\de_{j}\de_{j}^{*})^{k-1}),
    $$
    and
    \begin{align*}
	Tr( E(\cell_{j}K_{n}) (\de_{j}^{*}\de_{j})^{k}) & = Tr(
	|\de_{j}| E(\cell_{j}K_{n}) (\de_{j}^{*}\de_{j})^{k-1}
	|\de_{j}|) \\
	& = Tr( \de_{j} E(\cell_{j} K_{n} ) (\de_{j}^{*}\de_{j})^{k-1}
	\de_{j}^{*} )\\
	& = Tr( E(B_1(\cell_{j-1}K_{n})) \de_{j}
	E(\cell_{j}K_{n}) \de_{j}^{*} (\de_{j}\de_{j}^{*})^{k-1}).
    \end{align*}
    Since $B_{1}(\cell_{j-1}K_{n})\setminus \O_{n} \subset 
    B_{1}(\cf_{\cg}(\cell_{j-1}K_{n}))$, we obtain
    \begin{align*}
	| Tr( &E(\cell_{j}K_{n}) (\de_{j}^{*}\de_{j})^{k}) - Tr(
	E(\O_n) (\de_{j}\de_{j}^{*})^{k}) |
	\\
	& \leq Tr( E( B_1(\cf_\cg(\cell_{j-1}K_n)) ) ) \norm{\de_{j}
	E(\cell_{j}K_{n}) \de_{j}^{*}(\de_{j}\de_{j}^{*})^{k-1} }\\
	& \leq \m\norm{\de_{j}\de_{j}^{*}}^{k} |\cell_{j-1}(K_{n})|
	\eps_{n},
    \end{align*}
    where $\eps_{n} =
    \displaystyle{\frac{|\cf_{\cg}(\cell_{j-1}K_{n})|}
    {|\cell_{j-1}K_{n}|}}$ and
    $|B_{1}(\cf(\cell_{j-1}K_{n}))| \leq \m |\cell_{j-1}K_{n}|
    \eps_{n}$.    Finally,
    \begin{align*}
	\Big| &\frac{ Tr( E(\cell_{j-1}K_{n}) (\de_{j}\de_{j}^{*})^{k})
	} {|\cell_{p}K_{n}|} - \frac{ Tr( E(\cell_{j}K_{n})
	(\de_{j}^{*}\de_{j})^{k}) } {|\cell_{p}K_{n}|} \Big| \\
	& \leq |\cell_{p}K_{n}|^{-1} \Bigl( |Tr( E(
	\cf_\cg(\cell_{j-1}K_{n}) ) ( \de_{j}\de_{j}^{*} )^{k})| \\
	&+ | Tr( E( \cell_{j}K_{n} ) ( \de_{j}^{*}\de_{j} )^{k}) - Tr(
	E( \O_n ) ( \de_{j}\de_{j}^{*} )^{k}) | \Bigr) \\
	& \leq \frac{|\cell_{j-1}K_{n}|}{|\cell_{p}K_{n}|}(\m+1)
	\norm{\de_{j}\de_{j}^{*}}^{k} \eps_{n}\to 0, \text{ as }
	n\to\infty,
    \end{align*}
    by Lemma \ref{lem:diffNorm}.
\end{proof}

\begin{cor}\label{NSequality}
    For any continuous bounded function $f:[0,\infty)\to\bc$ vanishing
    at zero, one has $\Tr_{j-1}(f(\de_{j}\de_{j}^{*})) =
    \Tr_{j}(f(\de_{j}^{*}\de_{j}))$.  In particular
    $$
    \Tr_{j-1}(e^{-t\de_{j}\de_{j}^{*}}) - \Tr_{j-1}(I)
    = \Tr_{j}(e^{-t\de_{j}^{*}\de_{j}}) - \Tr_{j}(I).
    $$
    Therefore
    \begin{align*}
	\b_{j-1}^{+}(M) - \Tr_{j-1}(I)&= \b_{j}^{-}(M)- \Tr_{j}(I),
	\qquad \a_{j-1}^{+}(M)=\a_{j}^{-}(M),\\
	\b_{j-1}^{+}(M,\partial M) - \Tr_{j-1}(I)& =
	\b_{j}^{-}(M,\partial M)- \Tr_{j}(I), \qquad
	\a_{j-1}^{+}(M,\partial M)=\a_{j}^{-}(M,\partial M).
    \end{align*}
\end{cor}

\begin{proof}
    The proof for $\a(\D_{j\pm})$ and $\b(\D_{j\pm})$ follows directly
    by the previous results. All the arguments above may be rephrased
    for the relative invariants, giving the corresponding equality.
\end{proof}

\begin{rem}
    Let us recall that in \cite{LoLu} Novikov-Shubin numbers have been
    associated to the boundary operator $\partial$, namely depend on
    an index varying from 1 to the dimension $p$ of the complex.  As a
    consequence of Corollary \ref{NSequality}, there are only $p$ 
    independent Novikov-Shubin numbers in our framework too.
    
    Concerning $L^{2}$-Betti numbers, they have been defined in
    \cite{LoLu} as $\G$-dimensions of the $L^{2}$-homology, hence
    coincide with the trace of the kernel of the full Laplacians.  The
    relations proved in this section show that the $\b_{j}^{\pm}$'s
    are completely determined by the $\b_{j}$'s.  Moreover such
    relations imply a further identity which is the basis of the
    following Theorem.
\end{rem}

\begin{thm}[Euler-Poincar\'e characteristic]\label{EulerPoincare}
    Let $M$ be a $p$-dimensional self-similar CW-complex, with
    exhaustion $\set{K_{n}}$.  Then
    $$
    \chi^{\cg}(M):=\sum_{j=0}^{p}(-1)^{j}\b_{j}(M) =
    \lim_{n}\frac{\chi(K_{n})}{|\cell_{p}(K_{n})|}.
    $$
\end{thm}
\begin{proof}
    Let us first observe that, by using Theorem \ref{thm:BettiAndNS} 
    and Corollary \ref{NSequality}, we get
    \begin{align*}
	\chi^{\cg}(M)
	&=\b_{0}^{+}(M)+\sum_{j=1}^{p-1}(-1)^{j}
	(\beta_{j}^{+}(M) + \b_{j}^{-}(M) -Tr^{\cg}_{j}(I))+ 
	(-1)^{p}\b^{-}_{p}(M)\\
	&=\sum_{j=0}^{p-1}(-1)^{j} \b_{j}^{+}(M) + (-1)^{p}\Tr_{p}(I)
	+ \sum_{j=1}^{p}  (-1)^{j} \b_{j-1}^{+}(M) 
	- \sum_{j=1}^{p}  (-1)^{j} \Tr_{j-1}(I)	\\
	&= \sum_{j=0}^{p}(-1)^{j}\Tr_{j}(I).	
    \end{align*}
    On the other hand, 
    \begin{equation*}
	\sum_{j=0}^{p}(-1)^{j}\Tr_{j}(I))
	=\lim_{n}\sum_{j=0}^{p}(-1)^{j}
	\frac{|\cell_{j}(K_{n})|}{|\cell_{p}(K_{n})|}
	=\lim_{n}\frac{\chi(K_{n})}{|\cell_{p}(K_{n})|}.
    \end{equation*}
    The thesis follows.
\end{proof}
    
 \begin{ex}\label{EPexamples}
     $(i)$ For the Gasket graph of figure \ref{fig:Gasket} we get 
     $|VK_{n}| = \frac12 3^{n}+\frac32$ and $|EK_{n}| =3^{n}$, so 
     that $\chi^{\cg}(X)=-\frac12$.\\
     $(ii)$ For the Vicsek graph of figure \ref{fig:Vicsek} we get
     $|VK_{n}| = 3\cdot 5^{n}+1$ and $|EK_{n}| =4\cdot 5^{n}$, so that
     $\chi^{\cg}(X)=-\frac14$.\\
     $(iii)$ For the Lindstrom graph of figure \ref{fig:Lindstrom} we
     get $|VK_{n}| = 4\cdot 7^{n}+2$ and $|EK_{n}| =6\cdot 7^{n}$, so
     that $\chi^{\cg}(X)=-\frac13$.\\
 \end{ex}

\section{Prefractals as CW-complexes}\label{sec:prefractals}

We say that a $j$-dimensional polyhedron in some Euclidean space $\br^{m}$
is strictly convex if it is convex and any $(j-1)$-hyperplane contains
at most one of its $(j-1)$-dimensional faces.

\begin{defn}
    A polyhedral complex is a regular CW complex whose topology is
    that of a closed subset of some Euclidean space $\br^{m}$ and whose
    $j$-cells are flat strictly convex $j$-polyhedra in $\br^{m}$.  
\end{defn}

The following Proposition motivates the name of the boundary
subcomplex.

\begin{prop}
    Let $M$ be a $p$-dimensional polyhedral complex in $\br^{p}$.
    Then the boundary subcomplex $\partial M$ gives a CW structure on
    the boundary of $M$, seen as a subspace of $\br^{p}$.
\end{prop}

\begin{proof}
    Clearly, given a $(p-1)$-cell $\s$ of $M$ and a point $x\in\s$,
    we may find $\eps>0$ such that the ball $B(x,\eps)$ in $\br^{p}$
    is contained in $\s\cup\t_{1}\cup\t_{2}$, if $\s$ is contained in
    the two distinct $p$-cells $\t_{1},\t_{2}$, and is not contained
    in $M$ (indeed half of it is in the complement of $M$ w.r.t.
    $\br^{p}$), if $\s$ is contained in only one $p$-cell $\t$.  This
    proves the thesis.
\end{proof}

\begin{prop}\label{DistOnPoly}
    Let $M$ be a $p$-dimensional polyhedral complex, $j=1,\ldots,p$.
    If $\sigma\, ,\sigma^{\prime}$ are distinct cells in $\cell_{j-1}
    (M)$, there exists at most one polyhedron $\tau\in\cell_{j} (M)$
    such that $(\sigma,\partial_j \t)(\sigma',\partial_j \t)\ne0$.  \\
    If $\sigma\, ,\sigma^{\prime}\in \cell_{j} (M)$, there
    exists at most one polyhedron $\r\in\cell_{j-1} (M)$ such that
    $(\partial_j\sigma, \r)(\partial_j\sigma',\r)\ne0$.
\end{prop}
\begin{proof}
    Let $\sigma\neq\sigma^{\prime}\in \cell_{j-1} (M)$ and $C$ denote
    the convex hull of $\sigma\cup\sigma^{\prime}$.  If $C$ has
    dimension $j-1$, a $\t$ as above would have two faces in the same
    $(j-1)$-plane, against the strict convexity.

    \noindent If $C$ has dimension greater or equal to $j$, two
    different $\tau\, ,\tau^{\prime}\in \cell_{j} (M)$ containing both
    $\sigma$ and $\sigma^{\prime}$ in their boundaries would contain
    $C$ also, implying $\tau = \tau^{\prime}$.

    \noindent Finally, if $\sigma\, ,\sigma^{\prime}\in \cell_{j}
    (M)$, and there were $\r,\r'\in \cell_{j-1}(M)$ as above, they
    would belong to the $(j-1)$-plane separating $\s$ and $\s'$,
    namely $\s$ and $\s'$ would have two faces in the same
    $(j-1)$-plane, against strict convexity.
\end{proof}

Our main examples of self-similar CW-complexes will be a special
class of polyhedral complexes, namely prefractal complexes.

Let us recall that a self-similar fractal $K$ in $\br^{p}$ is
determined by contraction similarities $w_{1},\ldots,w_{q}$ as the
unique (compact) solution of the fixed point equation $K=WK$, where
$W$ is a map on subsets defined as $W A= \displaystyle{
\bigcup_{j=1}^{q} w_{j}A}$.

The fractal $K$ satisfies the {\it open set condition} with open set
$U$ if
$$
w_{j}U\subset U,\quad w_{j}U\cap w_{i}U=\emptyset,\quad i\ne j.
$$

Assume now we are given a self-similar fractal in $\br^{p}$ determined
by similarities $w_{1},\ldots,w_{q}$, with the same similarity
parameter.  Assume Open Set Condition holds for a bounded open set
whose closure is a strictly convex $p$-dimensional polyhedron $\cp$. 
If $\s$ is a multiindex of length $n$, we set $w_{\s}:=
w_{\s_{n}}\circ\cdots\circ w_{\s_{1}}$.  If $\s$ is an infinite
multiindex, we denote by $\s|{n}$ its $n$-th truncation.  Assume that
$w_{\s}\cp\cap w_{\s'}\cp$ is a (facial) subpolyhedron of both
$w_{\s}\cp$ and $w_{\s'}\cp$, $|\s|=|\s'|$.  Finally we choose an
infinite multiindex $I$.  We construct a polyhedral CW-complex as
follows.  First set 
$$
K_{n}:= 
w_{I|{n}}^{-1}W^{n}\cp=\bigcup_{|\s|=n}w_{I|{n}}^{-1}w_{\s}\cp.
$$

\begin{lem}
    $K_{n}$ is a finite polyhedral complex satisfying the following
    properties:
    \begin{align}
	\forall j< p, \ \s\in\cell_{j}(K_{n}),\  &\exists
	\t\in\cell_{p}(K_{n}):\s\subset\dot\t,\label{propF}\\
	K_{n}&\subset K_{n+1}.\label{increasing}
    \end{align}
\end{lem}
\begin{proof}
    Observe that $w_{I|{n}}^{-1}w_{\s}\cp \cap
    w_{I|{n}}^{-1}w_{\s'}\cp = w_{I|{n}}^{-1}(w_{\s}\cp \cap
    w_{\s'}\cp)$, hence is a common facial subpolyhedron, namely
    $K_{n}$ has a natural structure of polyhedral complex.  Property
    (\ref{propF}) is obvious.  Let us now prove (\ref{increasing}).
    \begin{equation*}
    K_{n+1} = w_{I|{n}}^{-1}w_{I_{n+1}}^{-1}\left(\bigcup_{j=1}^{q}
    w_{j}\right)W^{n}\cp \supset w_{I|{n}}^{-1}W^{n}\cp = K_{n}.
    \end{equation*}
\end{proof}

\begin{cor}
    $M=\cup_{n\in\bn} K_{n}$ is a polyhedral complex satisfying
    property (\ref{propF}). $\set{K_{n}}$ is an exhaustion for $M$.
    $M$ is called a prefractal complex.
\end{cor}
\begin{proof}
    Obvious.
\end{proof}

\begin{prop}
    A prefractal complex is a regular bounded CW-complex.
\end{prop}

\begin{proof}
    Regularity is obvious by construction.

    Let us estimate $V_{j}^{+}$.  Any $\s\in\cell_{j}$ is contained in
    some copy of the fundamental polyhedron $\cp$, therefore any
    $\t\in\ce_{j+1}$ such that $\dot{\t}\supset\s$ will be contained
    in the same copy of $\cp$ or in some neighboring copy.  As a
    consequence, we may estimate the number of such $\t$'s with the
    product of $|\cell_{j+1}(\cp)|$ times the maximum number of
    disjoint copies of $\cp$ having a $j$-cell in common.  Since such
    copies are contained in a ball of radius $\diam(\cp)$, their
    number may be estimated e.g. by the ratio of the volume of the
    ball of radius $\diam(\cp)$ and the volume of $\cp$.

    As for $V_{j}^{-}$, again any $\s\in\cell_{j}$ is contained in
    some copy of the fundamental polyhedron $\cp$, therefore any
    $\r\in\ce_{j-1}$ such that $\r\subset\dot{\s}$ will be contained
    in the same copy of $\cp$. The number of such $\r$'s is majorised
    by $|\cell_{j-1}(\cp)|$.
\end{proof}

In order to show that $M$ is a self-similar complex, we shall prove
that $K_{n}$ is a regular exhaustion satisfying  Definition
\ref{def:Quasiperiodic}.

\begin{lem}\label{lem:propF}
    Assume $K$ is a $p$-dimensional polyhedral complex in $\br^{p}$
    satisfying property (\ref{propF}), and we have
    $\s_{i}\in\cell_{i}(K)$, $i=j_{0},\ldots,j_{1}$, $j_{1}<p$, with
    $\s_{i}\subset\dot\s_{i+1}$.  Then there exist
    $\s'_{i}\in\cell_{i}(K)$, $i=j_{0},\ldots,j_{1}+1$, s.t.
    
    \itm{\a} $\s'_{i}\subset\dot\s'_{i+1}$, $i=j_{0},\ldots,j_{1}$,
    
    \itm{\b} $\s_{i-1}\subset\dot\s'_{i}$, $i=j_{0}+1,\ldots,j_{1}+1$,
    
    \itm{\g} $\text{ph}(\s_{i}'\cup\s_{j_{1}}) = \s'_{j_{1}+1}$,
    $i=j_{0},\ldots,j_{1}$, where $\text{ph}(E)$ denotes the
    polyhedral hull, namely the minimal polyhedron in $K$ (if it
    exists) containing $E\subset K$.
\end{lem}
\begin{proof}
    The proof will be done by descending induction on $i$, starting
    from $j_{1}+1$. \\
    Take any $\s'_{j_{1}+1}\supset \s_{j_{1}}$, which exists, because
    of (\ref{propF}).  Assume the statement for $k+1$.  If $k>j_{0}$,
    as $\s_{k-1} \subset \dot\s_{k}$, $\s_{k} \subset \dot\s_{k+1}'$,
    there is (one and only one) $\s'_{k}\subset\dot\s'_{k+1}$, s.t.
    $\s_{k-1}\subset\dot\s'_{k}$, and $\s_{k}'\neq \s_{k}$, by
    regularity of the CW-complex.  If $k=j_{0}$, simply take $\s'_{k}$
    as any $k$-dimensional face in the boundary of $\s'_{k+1}$,
    distinct from $\s_{k}$.  Let us observe that, for $i=j_{0}$,
    property $(\b)$ is empty.  This gives $(\a)$ and $(\b)$.  Since
    $\s'_{k+1}$ is strictly convex, and $\s_{k}$, $\s'_{k}$ are two
    distinct $k$-dimensional faces, we have $\text{ph}(\s'_{k}\cup
    \s_{k}) = \s'_{k+1}$.  Then $\text{ph}(\s_{j_{1}}\cup\s'_{k}) =
    \text{ph}(\s_{j_{1}} \cup \s_{k} \cup\s'_{k}) =
    \text{ph}(\s_{j_{1}}\cup\s'_{k+1}) = \s'_{j_{1}+1}$, which is
    property $(\g)$.
\end{proof}

\begin{thm}
    Any prefractal complex $M$ is a self-similar polyhedral complex.
\end{thm}
\begin{proof}
    We already proved regularity and boundedness.  Let us now observe
    that
    \begin{equation*}
	K_{n+1} = \bigcup_{j=1}^{q} w_{I|{n+1}}^{-1}w_{j}W^{n}\cp
	= \bigcup_{\ell=1}^{q} \g^{n}_{\ell}K_{n},
    \end{equation*}
    where $\g^{n}_{\ell}=w_{I|{n+1}}^{-1}w_{\ell}w_{I|{n}}$ are local
    isomorphisms.  Moreover, by OSC,
    $$
    \g^{n}_{k} K_{n} = w_{I|{n+1}}^{-1}w_{k}W^{n}\cp \subset
    w_{I|{n+1}}^{-1}w_{k}\cp.
    $$
    Therefore, $\g^{n}_{k}K_{n} \cap \g^{n}_{\ell} K_{n} \subset
    w_{I|{n+1}}^{-1} (w_{k} \cp \cap w_{\ell}\cp)$, which is contained
    in some affine hyperplane $\pi$; hence $\g^{n}_{k}K_{n} \cap
    \g^{n}_{\ell} K_{n}$ has no $p$-dimensional polyhedra.  Let
    $\s_{j_{0}}$ be a $j_{0}$-dimensional polyhedron in
    $\g^{n}_{k}K_{n} \cap \g^{n}_{\ell} K_{n}$, and $\s_{i}$,
    $i=j_{0},\ldots,j_{1}$, be a maximal family of polyhedra in
    $\g^{n}_{k}K_{n} \cap \pi$ s.t. $\s_{i}\subset \dot\s_{i+1}$,
    $i=j_{0},\ldots,j_{1}-1$.  Now apply the lemma with
    $K=\g^{n}_{k}K_{n}$, and get
    $\s'_{i}\in\cell_{i}(\g^{n}_{k}K_{n})$, $i=j_{0},\ldots,j_{1}+1$,
    with $\text{ph}(\s'_{i}\cup \s_{j_{1}}) = \s'_{j_{1}+1}$.  Then
    $\s'_{j_{0}}\not\subset\pi$, otherwise $\pi \supset \text{ph}
    (\s'_{j_{0}}\cup \s_{j_{1}}) = \s'_{j_{1}+1}$, against the
    maximality.  Therefore $\s'_{j_{0}}\not\in \cell_{j_{0}}(
    \g^{n}_{\ell} K_{n})$.  Moreover, since $\s_{j_{0}}$,
    $\s'_{j_{0}}\subset \dot\s'_{j_{0}+1}$, we have $d_{+}(\s_{j_{0}},
    \s'_{j_{0}})=1$, namely $\s_{j_{0}}\in
    \cf(\cell_{j_{0}}(\g^{n}_{\ell}K_{n}))$.  Similarly, one proves
    $\s_{j_{0}}\in \cf(\cell_{j_{0}}(\g^{n}_{k}K_{n}))$.  Therefore,
    $\g^{n}_{k}K_{n} \cap \g^{n}_{\ell} K_{n} = \cf(\g^{n}_{k}K_{n})
    \cap \cf(\g^{n}_{\ell} K_{n})$.
\end{proof}

\begin{rem}
    The construction above can be easily generalised to the
    case of translation limit fractals \cite{GuIs8, GuIs9, GuIs10,
    GuIs11, GuIs13}.
\end{rem}

We now study some properties of polyhedral complexes, which are valid
in particular for the prefractal complexes.

\begin{defn}
    We say that $\Delta_{j\pm}$ is a {\it graph-like Laplacian} if
    there exists a suitable orientation of $M$ such that the
    off-diagonal entries of the matrix associated with
    $\Delta_{j\pm}$, in the corresponding orthonormal basis, belong to
    $\{0,-1\}$.
\end{defn}

\begin{thm}\label{Thm:10.2}
    Assume $M$ is a $p$-dimensional polyhedral complex in $\br^{m}$,
    $m\geq p$.  Then $\Delta_{j+}$ is graph-like if and only if $j=0$.
\end{thm}
\begin{proof}
    If $\s,\s'\in\cell_{0}(M)$, then
    $(\s,\D_{0+}\s')=\sum_{\a}[\t^{1}_{\a},\s][\t^{1}_{\a},\s']$.  By
    the Proposition \ref{DistOnPoly}, the sum consists of at most one
    non-vanishing summand, corresponding to some 1-cell $\t$.  By the
    choice we made for the orientation on the 0-cells, the sum of
    $[\t^{1}_{\a},\s]$ and $[\t^{1}_{\a},\s']$ is 0, hence the product
    is $-1$.

    Let $j>0$ and choose $\tau\in \cell_{j+1} (M)$.  Since $j+1 \ge
    2$, $\tau$ has at least three distinct faces $\sigma_{1}\,
    ,\sigma_{2}\, ,\sigma_{3}\in \cell_{j} (M)$.  Setting
    $[\t,\s_{j}]=\l_{j}=\pm1$, we obtain, for $j\ne k$,
    $\bigl(\sigma_{i} ,\Delta_{j+}\sigma_{k} \bigr)
    =\lambda_{i}\cdot\lambda_{k}$. Therefore the product
    $$\bigl(\sigma_{1} ,\Delta_{j+}\sigma_{2} \bigr)\cdot
    \bigl(\sigma_{2} ,\Delta_{j+}\sigma_{3} \bigr)\cdot
    \bigl(\sigma_{3} ,\Delta_{j+}\sigma_{1} \bigr)=1$$
    so that the three off diagonal matrix elements will never be all
    equal to $-1$.
\end{proof}

\begin{lem}\label{DeltaMinusNograph}
    Assume $\D_{(j+1)-}$ is not diagonal. Then $\D_{j-}$ is not
    graph-like.
\end{lem}

\begin{proof}
    By assumption, there exist $\t_1,\t_2\in\cell_{j+1}(M)$ such that
    $(\partial_{j+1}\t_1,\partial_{j+1}\t_2)= (\t_{1},
    \D_{(j+1)-}\t_{2})\ne0$, namely there exist
    $\s_3\subset\dot{\t}_1\cap\dot{\t}_2$.  If $\r\subset\dot{\s}_3$,
    $\partial_j(\partial_{j+1}\t_i)=0$ implies the existence of
    $\s_i\subset\dot{\t}_i$ for which $\r\in\dot{\s}_i$.  Setting
    $\l_i=[\s_i,\r]$, $i=1,2,3$, $(\s_i,\D_{j-}\s_k)=\l_i\l_k$, the
    proof goes on as in Theorem \ref{Thm:10.2}.
\end{proof}

In order to prove a general result on the possibility of
$\Delta_{j-}$ to be graph-like, we shall exclude some trivial
cases.

\begin{defn}\label{dfn:irreducible}
    We say that an operator $A$ acting on $\ell^2 (\cell_j(M))$ is
    irreducible if the only self-adjoint idempotent multiplication
    operators commuting with $A$ are $0$ and $I$.  We say that a
    (connected) polyhedral complex is $q$-irreducible if $q$ is the
    maximum number such that $\D_{j+}$ is irreducible for any $j<q$,
    and $\D_{j-}$ is irreducible for any $j\leq q$.
\end{defn}

\begin{thm}\label{lem:notgraph}
    Let $M$ be a $q$-irreducible $p$-dimensional polyhedral complex in
    $\br^p$, $j\leq q$.  Then $\D_{j-}$ graph-like implies $j=q$.  If
    $q=p$ the above implication is indeed an equivalence.
\end{thm}
\begin{proof}
    If $j<q$, $\D_{j-}$ is not graph-like by Lemma
    \ref{DeltaMinusNograph}.  If $q=p$ we may choose the orientation
    for $p$-cells according to a given orientation of $\br^p$.  Then,
    if $\tau,\tau'\in\cell_{p}(M)$ have a face in common, i.e.
    $(\tau,\Delta_{p-}\tau')\ne 0$, such face receives opposite
    orientations from the embeddings in $\dot\tau$, resp.  $\dot
    \tau'$, i.e. $(\tau,\Delta_{p-}\tau') = (\partial_p\tau,
    \partial_p\tau') = -1$.
\end{proof}

\begin{defn}\label{dualgraph}
    If $M$ is a $p$-irreducible polyhedral complex in $\br^{p}$, by
    the previous Theorem, a graph $G$ is associated to $\ov{\D}_{p-}$,
    and is constructed as follows.  The set of vertices of $G$ is
    $\cell_{p}(M)$, while $\s,\s'\in\cell_{p}(M)$ are adjacent {\it
    iff} there is $\r\in\cell_{p-1}(M)$ such that
    $\r\subset\dot{\s},\dot{\s}'$.  We call $G$ the dual graph of $M$.
\end{defn}

\section{Computation of the Novikov-Shubin numbers for fractal graphs}

Let us observe that a $1$-dimensional regular CW-complex is the same
as a simple graph, 
and boundedness means bounded degree.  Recall that a simple graph  $G=(VG,EG)$ is a collection $VG$ of objects,
 called  {\it vertices}, and a collection $EG$ of unordered pairs of
 distinct vertices, called {\it edges}. We call a self-similar
$1$-dimensional CW-complex simply a self-similar graph.

The results in this section will allow us to calculate $\a_0=\a_{1}$
of some self-similar graphs, and also $\a_p$ of a $p$-irreducible
prefractal complex in $\br^p$, in the sense of definition
\ref{dfn:irreducible}.

In the rest of this section, $G$ is a countably infinite graph with
bounded degree.  We denote by $\D$ the Laplacian on $0$-cells
(points), hence $\D=C-A$, where $C$ is a diagonal operator, with
$C(x,x)=c(x)$ the number of edges starting from the point $x$, and $A$
is the adjacency matrix.  Let $P$ be the transition operator, i.e.
$p(x,y)$ is the transition probability from $x$ to $y$ of the simple
random walk on $G$.  Let $\ca_{0}$ be a C$^{*}$-algebra of operators,
acting on $\ell^{2}(G)=\ell^{2}(\cell_{0}(G))$, which contains $\D$,
and possesses a finite trace $\t_0$.



We can also consider the Hilbert space $\ell^{2}(G,c)$ with scalar
product $(v,w)_{c}=\sum_{x\in G}c(x)\ov{v}_{x}w_{x}$.  On this space
the transition operator $P$ is selfadjoint, and the Laplace operator
is defined as $\D_{c}=I-P$.  Since $G$ has bounded degree, $C$ is
bounded from above by a multiple of the identity, and, since $G$ is
connected, it is bounded from below by the identity.  Also, the two
spaces $\ell^{2}(G)$ and $\ell^{2}(G,c)$ coincide as topological
vector spaces, with the obvious identification.  With this
identification we have $(u, \D v)=(u,\D_{c}v)_{c}$.  We may also
identify operators on the two spaces, hence the C$^{*}$-algebras
$\ca_{0}$, resp.  $\ca_{0,c}$, acting on $\ell^{2}(G)$, resp.
$\ell^{2}(G,c)$, can be identified as topological algebras.  We use
this identification to carry the trace $\t_0$ onto $\ca_{0,c}$.

\begin{thm}\label{Delta-DeltaC}
    Let $\m$ be the maximal degree of $G$. Then
    $$
    \t_0\left(\frac{1}{1+\m t\D_{c}}\right) \leq
    \t_0\left(\frac{1}{1+t\D}\right) \leq
    \t_0\left(\frac{1}{1+t\D_{c}}\right) \qquad t\geq 0.
    $$
\end{thm}
\begin{proof}
    Let us consider the positive self-adjoint operator $Q=C^{-1/2}\D
    C^{-1/2}$ on $\ell^{2}(G)$.  Since $C$ and $\D$ belong to
    $\ca_{0}$, $Q\in\ca_{0}$ too.  Since $A=CP$, we have $Q=C^{1/2}
    (I-P)C^{-1/2}$. Clearly $Q\leq Q^{1/2}CQ^{1/2}\leq\m Q$, hence, by
    operator monotonicity,
    $$
    \t_0\left(\frac{1}{1+\m tQ}\right)
    \leq \t_0\left(\frac{1}{1+tQ^{1/2}CQ^{1/2}}\right)
    \leq \t_0\left(\frac{1}{1+tQ}\right)
    \qquad t\geq 0.
    $$
    Now observe that, for sufficiently small $t$,
    \begin{align*}
    \t_0\left(\frac{1}{1+tQ}\right)
    &=\sum_{n=0}^{\infty}(-t)^{n}\t_0(Q^{n})\\
    &=\sum_{n=0}^{\infty}(-t)^{n}\t_0(\D_{c}^{n})
    =\t_0\left(\frac{1}{1+t\D_{c}}\right).
    \end{align*}
    Since both the left and the right hand side are analitic functions
    for $t\geq0$, they coincide for any $t\geq 0$.  Analogously,
    \begin{align*}
    \t_0\left(\frac{1}{1+tQ^{1/2}CQ^{1/2}}\right)
    &=\sum_{n=0}^{\infty}(-t)^{n}\t_0((Q^{1/2}CQ^{1/2})^{n})\\
    &=\sum_{n=0}^{\infty}(-t)^{n}\t_0(\D^{n})
    =\t_0\left(\frac{1}{1+t\D}\right).
    \end{align*}
    The result follows.
\end{proof}

In order to prove the main result of this section, we need a Tauberian theorem.  It is a quite simple modification
of a theorem of de Haan and Stadtm\"{u}ller, cf.  \cite{BGT} thm.
2.10.2, and, on the same book, also thm.  1.7.6, by Karamata, showing
that the bound $\a<1$ below is a natural one.

\begin{defn}
    \itm{i} Let us denote by $OR(1)$ the space of positive, non
    increasing functions $f$ on $[0,\infty)$ such that $\exists
    T>0,c>0,\a\in(0,1)$ such that
    \begin{equation}\label{OR1}
    \frac{f(\l t)}{f(t)}\leq c\l^{-\a},\qquad \forall \l>0,\ t\geq T.
    \end{equation}
    
    \itm{ii} If $f,g$ are functions on $[a,+\infty)$ we write
    $f\asymp g$, $t\to\infty$ if $\exists
    T\geq a,k>1$ such that
    \begin{equation}\label{pesci}
    k^{-1}g(t)\leq f(t)\leq k g(t),\qquad t\geq T.
    \end{equation}
\end{defn}

\begin{rem}\label{rem:allT}
    If the functions are bounded and defined on $[0,\infty)$, we may
    equivalently assume $T=0$ both in eq.  (\ref{OR1}) and
    (\ref{pesci}), possibly changing the constants $c$ and $k$.
\end{rem}

Let us now denote by $\hat{f}$ the Laplace transform of $f$,
$$
\hat{f}(t)=t\int_{0}^{\infty}e^{-ts}f(s)\ ds.
$$

\begin{lem}\label{Tauberian}
    Let $f$ be a positive bounded function. Then
    $f\in OR(1)$ {\it {iff}} $\hat{f}(1/\cdot)\in OR(1)$. In this case
    $$
    f\asymp \hat{f}(1/\cdot),\qquad t\to\infty.
    $$
\end{lem}
\begin{proof}
    Let us notice that since $f$ is bounded also $\hat{f}$ is
    bounded, hence, according to Remark \ref{rem:allT}, the properties
    above should hold for all $t\geq0$.
    Now observe that
    \begin{equation*}\label{eq:Laplace}
    \hat{f}(s)=\int_{0}^{\infty}e^{-y}f(y/s)\ dy,
    \end{equation*}
    hence
    \begin{equation}\label{eq:Laplace2}
    \hat{f}(1/s)=\int_{0}^{\infty}e^{-y}f(sy)\ dy
    \geq f(xs)\int_{0}^{x}e^{-y}\ dy,
    \end{equation}
    giving
    \begin{equation}\label{ineq1}
    f(t)\leq\frac{e^{x}}{e^{x}-1}\hat{f}(x/t), \forall x,t.
    \end{equation}
    Now assume $f\in OR(1)$.
    \\
    By (\ref{eq:Laplace2}) we get
    $\displaystyle{\frac{\hat{f}(1/t)}{f(t)}
    =\int_{0}^{\infty}e^{-\l}\frac{f(\l t)}{f(t)}\ d\l}$. Therefore,
    splitting the domain of integration and using property $OR(1)$, we
    get
    $$
    \frac{\hat{f}(1/t)}{f(t)}\leq c\int_{0}^{1}\l^{-\a}e^{-\l}\ d\l +
    \int_{1}^{\infty}e^{-\l}\ d\l\leq c\Gamma(1-\a)+1.
    $$
    This, together with (\ref{ineq1}) for $x=1$, implies $f\asymp
    \hat{f}(1/\cdot)$, $t\to\infty$.
    \\
    Moreover,
    \begin{equation*}\label{isOR1}
    \frac{\hat{f}(1/(\l t))}{\hat{f}(1/t)}
    \leq c' \frac{f(\l t)}{\hat{f}(1/t)}\leq c'' \frac{f(\l t)}{f(t)}
    \leq c'''\l^{-\a},
    \end{equation*}
    showing that $\hat{f}(1/\cdot)\in OR(1)$.
    \\
    Now assume $\hat{f}(1/\cdot)\in OR(1)$. Then, by
    (\ref{ineq1}) with $x=1$,
    \begin{align*}
    \hat{f}(s)&=\int_{0}^{a}e^{-y}f(y/s)\ dy
    + \int_{a}^{\infty}e^{-y}f(y/s)\ dy\\
    &\leq \frac{e}{e-1} \int_{0}^{a}e^{-y}\hat{f}(s/y)\ dy
    + f(a/s)\int_{a}^{\infty}e^{-y}\ dy\\
    &\leq c \frac{e}{e-1} \hat{f}(s)\int_{0}^{a}y^{-\a}e^{-y}\ dy
    + f(a/s).
    \end{align*}
    Choosing $a$ sufficiently small, we get $c \frac{e}{e-1}
    \int_{0}^{a}y^{-\a}e^{-y}\ dy \leq1/2$, hence
    $\hat{f}(s)\leq2f(a/s)$. Now, using (\ref{ineq1}) with
    $x=a$, we get
    $$
    \frac{f(\l t)}{f(t)}\leq\frac{2e^{a}}{e^{a}-1}
    \frac{\hat{f}\left(\frac{a}{\l
    t}\right)}{\hat{f}\left(\frac{a}{t}\right)} \leq
    \frac{2c'e^{a}}{e^{a}-1}\l^{-\a},
    $$
    namely $f\in OR(1)$. The proof is complete.
\end{proof}

\begin{cor}\label{exp-res}
    Let $(\ca,\t)$ be a C$^{*}$-algebra with a finite trace, $A$ be a
    positive element of $\ca$.  Then $\t(e^{-tA})\in OR(1)$ {\it
    {iff}} $\t((1+tA)^{-1})\in OR(1)$.  In this case
    $\t(e^{-tA})\asymp \t((1+tA)^{-1})$, $t\to\infty$.
\end{cor}
\begin{proof}
    Let us consider, in the von Neumann algebra of the GNS
    representation of $\t$, the function $N(t)=\t(e_{[0,t)}(A))$. 
    Then $f(t) :=\t(e^{-tA}) = \int_{0}^{\infty}e^{-ts}\ dN(s)$, hence
    its Laplace transform is
    $$
    \hat{f}(x) = x\int_{0}^{\infty}dN(s)\int_{0}^{\infty}dt\
    e^{-t(s+x)} =\int_{0}^{\infty}\frac{x}{s+x}\ dN(s)
    =\t\left(\frac{1}{1+x^{-1}A}\right).
    $$
    The result now follows from the Lemma above.
\end{proof}

Now we come back to the Laplacians on $G$.

\begin{cor}\label{kernelasym}
    Let $G$ be a countably infinite connected graph with bounded
    degree.  Let $\D$, resp.  $\D_{c}$, be the homological, resp. 
    probabilistic Laplacian.  Let $\ca_{0}$ be a C$^{*}$-algebra of
    operators, acting on $\ell^{2}(G)=\ell^{2}(\cell_{0}(G))$, which
    contains $\D$ (hence $\D_c$), and possesses a finite trace $\t_0$. 
    Consider the functions $\th(t)=\t_0(e^{-t\D})$,
    $\th_{c}(t)=\t_0(e^{-t\D_{c}})$.  Then $\th\in OR(1)$ {\it {iff}}
    $\th_{c}\in OR(1)$.  In this case $\th\asymp \th_{c}$ for
    $t\to\infty$.
\end{cor}
\begin{proof}
    Assume $\th\in OR(1)$.  Then, by Corollary \ref{exp-res},
    $\th(t)\asymp \t_0((1+t\D)^{-1})$,
    $t\to\infty$ and $\t_0((1+t\D)^{-1})\in OR(1)$.  Therefore,
    recalling Theorem \ref{Delta-DeltaC}, and denoting by $\m$ the
    maximum degree of $G$, we get
    $$
    1\leq \frac{\t_0((1+t\D_{c})^{-1})}{\t_0((1+t\D)^{-1})}
    \leq \frac{\t_0((1+\m^{-1}t\D)^{-1})}{\t_0((1+t\D)^{-1})}
    \leq c\m^{\a},
    $$
    namely $\t_0((1+t\D_{c})^{-1})\asymp \t_0((1+t\D)^{-1})$.
    Therefore $\t_0((1+t\D_{c})^{-1})\in OR(1)$.  Applying
    Corollary \ref{exp-res} again we get $\th_{c}\in OR(1)$ and
    $\th\asymp \th_{c}$.  The converse implication is proved
    analogously.
\end{proof}

Now we relate the large $n$ asymptotics of the probability of
returning to a point in $n$ steps with the large time heat kernel asymptotics. 
Since for bipartite graphs the probability is zero for odd $n$, the
estimates are generally given in terms of the sum of the $n$-step plus
the $(n+1)$-step return probability.  In order to match the above
treatment we shall use a suitable mean for the return probability,
namely the trace of the $n$-th power of the transition operator $P$.

First we need the auxiliary function described in the following

\begin{lem}
    Let us denote by $\f_{\g}$, $\g>0$, the function
    $$\f_{\g}(x):=e^{x}x^{-\g}\int_{0}^{x}e^{-t}d(t^{\g}),\, x\geq0.$$
    Then $\f_{\g}$ extends to the entire function
    \begin{equation}\label{powerseries}
    \sum_{n=0}^{\infty}\frac{x^{n}}{n!}
    \begin{pmatrix}
        n+\g\\n
    \end{pmatrix}^{-1},\, x\in\bc.
    \end{equation}
\end{lem}
\begin{proof}
    Let us observe that the power series in (\ref{powerseries}) is an
    entire function $\f$ satisfying
    $$
    \begin{cases}
    \f'=(1-\frac{\g}{x})\f+\frac{\g}{x}\\
    \f(0)=1\, .
    \end{cases}
    $$
    It is easy to see that $\f_{\g}$ is the unique solution of the
    differential equation in $(0,+\infty)$ which tends to 1 when $x\to
    0^{+}$.
\end{proof}

\begin{thm}\label{main}
    Let $G$ be a countably infinite connected graph with bounded
    degree.  Let $\ca_{0}$ be a C$^{*}$-algebra of operators, acting
    on $\ell^{2}(G)=\ell^{2}(\cell_{0}(G))$, which contains the
    homological Laplacian $\D$ and possesses a finite trace $\t_0$. 
    If $\t_0(P^{n}+P^{n+1})\asymp n^{-\g}$, then $\th_{c}(t)\asymp
    t^{-\g}$, $t\to\infty$.
\end{thm}
\begin{proof}
    Let us observe that
    $$
    e^{-x}x^{\g}\f_{\g}(x)\to\g\G(\g),\quad x\to+\infty,
    $$
    and that
    $$
    \begin{pmatrix}
    n+\g\\n
    \end{pmatrix}
    n^{-\g}\to1,\quad n\to+\infty.
    $$
    On the one hand, we have
    \begin{align*}
	\th_{c}(t)&=e^{-t}\t_0(e^{tP})
	=e^{-t}\sum_{n=0}^{\infty}\frac{t^{n}}{n!}\t_0(P^{n}) \leq K
	e^{-t}\sum_{n=0}^{\infty}\frac{t^{n}}{n!}n^{-\g}\\
	&\leq K' e^{-t}\sum_{n=0}^{\infty}\frac{t^{n}}{n!}
	\begin{pmatrix}
	    n+\g\\n
	\end{pmatrix}^{-1}
	=K'e^{-t}\f_{\g}(t)\leq K''t^{-\g}.
    \end{align*}
    On the other hand, setting $\psi(t) = \sum_{n=0}^{\infty}
    \frac{t^{n}}{n!} \t_0(P^{n})$, we get $\th_{c}(t)=
    e^{-t}\psi(t)$, and $
    2\th_{c}(t)+\th'_{c}(t)=e^{-t}(\psi(t)+\psi'(t))$.  Let us note
    that
    $$
    \psi'(t)=\sum_{n=1}^{\infty}\frac{t^{n-1}}{(n-1)!}\t_0(P^{n})
    =\sum_{n=0}^{\infty}\frac{t^{n}}{n!}\t_0(P^{n+1}),
    $$
    therefore
    $$
    \psi(t)+\psi'(t) =
    \sum_{n=0}^{\infty}\frac{t^{n}}{n!}\t_0(P^{n}+P^{n+1}) \geq c
    \sum_{n=0}^{\infty}\frac{t^{n}}{n!}n^{-\g}\geq c' e^{t}t^{-\g}.
    $$
    Finally, since $\th_{c}'$ is negative,
    $$
    2\th_{c}\geq2\th_{c}(t)+\th_{c}'(t)=e^{-t}(\psi(t)+\psi'(t))\geq
    c' t^{-\g}.
    $$
    The thesis follows.
\end{proof}

\begin{cor}\label{cor:NSforGraphs}
    Let $G$ be a countably infinite connected graph with bounded
    degree. Let $\ca_{0}$ be a C$^{*}$-algebra  of operators, acting
    on $\ell^{2}(G)=\ell^{2}(\cell_{0}(G))$, which contains $\D$, and
    possesses a finite trace $\t_0$. Assume
    $\t_0(P^{n}+P^{n+1})\asymp n^{-\g}$, for $\g>0$.  
    
    \itm{i} If $\g\in(0,1)$, then the
    Novikov-Shubin number $\a(G)=2\g$.
    
    \itm{ii} If $\g>0$ and $G$ has constant degree, then $\a(G)=2\g$.
\end{cor}
\begin{proof}
    $(i)$ it follows from Corollary \ref{kernelasym} and Theorem
    \ref{main}.
    
    $(ii)$ it follows from Theorem \ref{main} and the observation that
    if the degree is constantly equal to $\m$, then $\t_0(e^{-t\D}) =
    \t_0(e^{-\m t\D_{c}})$.
\end{proof}

\begin{cor}\label{cor:HKestimates}
    Let $G$ be a countably infinite connected graph with bounded
    degree. Let $\ca_{0}$ be a C$^{*}$-algebra  of operators, acting
    on $\ell^{2}(G)=\ell^{2}(\cell_{0}(G))$, which contains $\D$, and
    possesses a finite trace $\t_0$. Denote by
    $p_{n}(x,y)$ the $(x,y)$-element of the matrix $P^{n}$, which is
    the probability that a simple random walk started at $x$ reaches
    $y$ in $n$ steps.  Assume that there are $\g\in(0,1)$,
    $c_{1},c_{2}>0$ such that, for all $x\in G$, $n\in\bn$,
    \begin{equation}\label{eq:HKestimates}
    \begin{split}
    p_{n}(x,x) & \leq c_{2}n^{-\g}\\
    p_{n}(x,x)+p_{n+1}(x,x) & \geq c_{1}n^{-\g}.
    \end{split}
    \end{equation}
    Then the Novikov-Shubin number $\a(G)=2\g$.
\end{cor}
\begin{proof}
    It is just a restatement of the previous Corollary.
\end{proof}

 As Novikov-Shubin numbers of covering manifolds are large scale
 invariants, one expects that graphs which are asymptotically close
 should have the same Novikov-Shubin number.  We show that this
 happens in case of roughly isometric graphs.
 
\begin{defn}
    Let $G_{1}$, $G_{2}$ be infinite graphs with bounded degree.  A
    map $\f:G_{1}\to G_{2}$ is called a {\it rough isometry} if there
    are $a,\, b,\, M>0$ s.t.
    
    \itm{i} $a^{-1} d_{1}(x,y) -b \leq d_{2}(\f(x),\f(y)) \leq a 
    d_{1}(x,y) +b$, for $x,y\in G_{1}$,
    
    \itm{ii} $d_{2}(\f(G_{1}),y) \leq M$, for $y\in G_{2}$.
    
    Then $G_{1}$ and $G_{2}$ are said to be rough isometric.
\end{defn}

 Observe that being rough isometric is an equivalence relation.

\begin{thm} \label{thm:invariance}
    Let $G_{1}$, $G_{2}$ be rough isometric, infinite graphs with
    bounded degree.  For $j=1,2$, let $\ca_{j}$ be a C$^{*}$-algebra
    of operators, acting on $\ell^{2}(G_j)=\ell^{2}(\cell_{0}(G_j))$,
    which contains the Laplace operator $\D_j$ of the graph $G_j$, and
    possesses a finite trace $\t_j$.  Assume $G_{1}$ satisfies
    (\ref{eq:HKestimates}), then $G_{2}$ does as well.  As a
    consequence, $\a_0(G_1)=\a_0(G_2)$.
\end{thm}
\begin{proof}
    It is a consequence of \cite{GrTe} Theorem 3.1, \cite{Barlow2}
    Lemma 1.1, and \cite{HaKu} Theorem 5.11.
\end{proof}

\section{Examples}

In this section, we compute the Novikov-Shubin numbers of some
explicit examples.

Our first class of examples is that of nested fractal graphs, for more
details on the construction see Section \ref{sec:prefractals} and
\cite{HaKu}.

Assume we are given a nested fractal $K$ in $\br^{p}$ determined by
similarities $w_{1},\ldots,w_{q}$, with the same similarity parameter,
and let $S$ be the Hausdorff dimension of $K$ in the resistance metric
\cite{Kiga}.  Let $M$ be the nested fractal graph based on $K$.

\begin{thm}\label{thm:HaKu}
    Let $K$, $S$, and $M$ be as above.  Then (\ref{eq:HKestimates}) 
    hold for $M$, with $\g=\frac{S}{S+1}\in(0,1)$. Therefore, 
    $\a_{0}(M) =2\g$.
\end{thm}
\begin{proof}
    The thesis follows from Corollary 4.13 in \cite{HaKu}, and
    Corollary \ref{cor:HKestimates} above.
\end{proof}

\begin{ex}
    Using the previous Theorem we can compute $\a_{0}$ for some
    self-similar graphs coming from fractal sets as in Section
    \ref{sec:prefractals}.  Moreover, since $\b_{0}=0$ by the estimate
    in Theorem \ref{main}, we get, by Corollary \ref{EulerPoincare},
    $\b_{1}=-\chi$, the latter being explicitely computed in Example
    \ref{EPexamples}.

    For the Gasket graph in figure \ref{fig:Gasket} we obtain 
     $\a_{0}= \frac{2\log3}{\log5}$, see \cite{Barlow1},
    $\b_{0}=0$, and $\b_{1}=\frac12$.

    For the Vicsek graph in figure \ref{fig:Vicsek} we obtain $\a_{0}=
    \frac{2\log 5}{\log 15}$, see \cite{KiLa},
    $\b_{0}=0$, and $\b_{1}=\frac14$.
    
    For the Lindstrom graph in figure \ref{fig:Lindstrom} we obtain
    $\a_{0}= \frac{2\log7}{\log 12.89027}$ computed numerically, see
    \cite{Kuma}, $\b_{0}=0$, and $\b_{1}=\frac13$.
\end{ex}

Our second class of examples is given by the following

\begin{prop}\label{Prop:DualGraph}
    Let $M$ be a $p$-irreducible prefractal complex in $\br^p$, let
    $G$ be the dual graph of $M$, as in Definition \ref{dualgraph},
    and assume that (\ref{eq:HKestimates}) hold on $G$.  Then the
    Novikov-Shubin number $\a_p(M,\de M) =2\g$.
\end{prop}
\begin{proof}
    It is a consequence of Corollary \ref{cor:HKestimates}.
\end{proof}

\begin{ex}
    Let us consider the 2-dimensional complex $M$ in figure 
    \ref{fig:Dodecagon}.
    \begin{figure}[ht]
    \centering
    \psfig{file=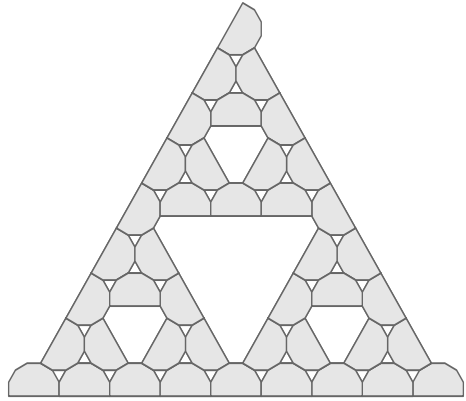,height=2.165in,width=2.5in}
    \caption{Dodecagon 2-complex}
    \label{fig:Dodecagon}
    \end{figure}
    We want to compute its second relative Novikov-Shubin number
    $\a_2(M,\partial M)$.  If we extend the definition of selfsimilar
    CW-complex as explained in Remark \ref{rem:rBdryIntersec}, it is
    easy to see that the complex in figure \ref{fig:Dodecagon} is
    self-similar and its dual graph, defined in Definition
    \ref{dualgraph}, coincides with the Gasket graph considered in
    figure \ref{fig:Gasket}.  Therefore, by Proposition
    \ref{Prop:DualGraph}, $\a_2(M,\partial M)= \frac {2\log3}
    {\log5}$.
\end{ex}

\begin{ex}
    The Carpet $2$-complex $M$ in figure \ref{fig:Carpet} is an
    example of a $2$-dimensional self-similar CW-complex.  Barlow
    \cite{Barlow1} associates to $M$ the graph $G_{1}$ in figure
    \ref{fig:BarlowCarpetGraph}, which, by \cite{Barlow1} Theorem 3.4,
    satisfies the estimates in Corollary \ref{cor:HKestimates}, with
    $\g=\frac{\log 8}{\log(8\r)}$, where $\r\in[\frac76,\frac32]$, while computer calculations suggest that $\r\cong 1.251$.

    The dual graph of $M$, as in Definition \ref{dualgraph}, is the
    graph $G_{2}$ in figure \ref{fig:dualCarpetGraph}, also associated
    to $M$ by Barlow and Bass in \cite{BaBa}.  It is easy to see that
    the graphs $G_{1}$ and $G_{2}$ are roughly-isometric.  Therefore,
    by Theorem \ref{thm:invariance}, $\a_{2}(M,\de M) =
    \a(\ov{\D}_{2-}) = 2\g$ (so that $\a_{2}(M,\de M) \in [1.67,
    1.87])$.
 
    \begin{figure}[ht]
    \centering
    \psfig{file=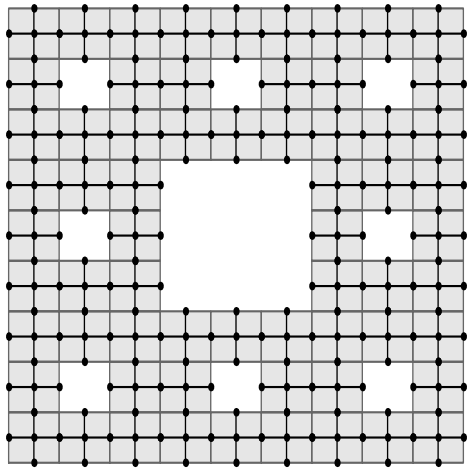,height=2.5in,width=2.5in}
    \caption{Carpet Graph (following Barlow)}
    \label{fig:BarlowCarpetGraph}
    \end{figure}

    \begin{figure}[ht]
    \centering
    \psfig{file=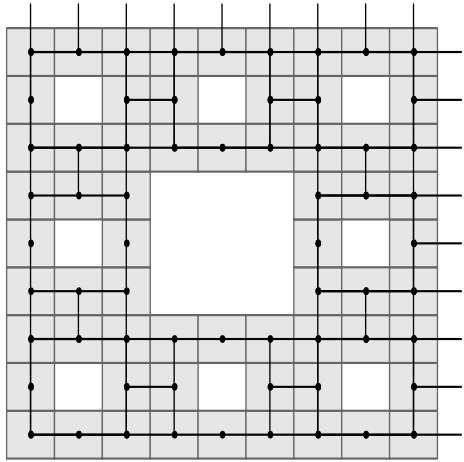,height=2.5in,width=2.5in}
    \caption{Dual Carpet Graph}
    \label{fig:dualCarpetGraph}
    \end{figure}
\end{ex}


\end{document}